\input amstex
\documentstyle{amsppt}
%
%
\nopagenumbers
\def\Img{\operatorname{Im}}
\def\id{\operatorname{id}}
\def\tr{\operatorname{tr}}
\accentedsymbol\tx{\tilde x}
\def\negskp{\hskip -2pt}
\def\compos{\,\raise 1pt\hbox{$\sssize\circ$} \,}
\pagewidth{360pt}
\pageheight{606pt}
\leftheadtext{Ruslan A. Sharipov}
\rightheadtext{Comparative analysis for pair of dynamical systems \dots}
\topmatter
\title Comparative analysis for pair of dynamical systems,
one of which is Lagrangian.
\endtitle
\author
R.~A.~Sharipov
\endauthor
\abstract
It is known that some equations of differential geometry are
derived from variational principle in form of Euler-Lagrange
equations. The equations of geodesic flow in Riemannian
geometry is an example. Conversely, having Lagrangian dynamical
system in a manifold, one can consider it as geometric
equipment of this manifold. Then properties of other dynamical
systems can be studied relatively as compared to this Lagrangian
one. This gives fruitful analogies for generalization. In present
paper theory of normal shift of hypersurfaces is generalized from
Riemannian geometry to the geometry determined by Lagrangian
dynamical system. Both weak and additional normality equations for
this case are derived.
\endabstract
\address Rabochaya street 5, 450003, Ufa, Russia
\endaddress
\email \vtop to 20pt{\hsize=280pt\noindent
R\_\hskip 1pt Sharipov\@ic.bashedu.ru\newline
r-sharipov\@mail.ru\vss}
\endemail
\urladdr
http:/\negskp/www.geocities.com/r-sharipov
\endurladdr
\endtopmatter
\loadbold
\TagsOnRight
\document
\head
1. Introduction.
\endhead
    Let $M$ be some smooth manifold and let $L$ be Lagrange
function of some Lagrangian dynamical system in $M$. This is
smooth scalar function depending on point $p$ of $M$ and on
vector $\bold v$ at this point. In other words, $L$ is a
function of point $q=(p,\bold v)$ of tangent bundle $TM$.
We treat $L$ as basic equipment of manifold $M$ (like metric
tensor in Riemannian geometry or simplectic structure in
simplectic geometry). Therefore we assume $L$ to be satisfying
some special requirements. Let $x^1,\,\ldots,\,x^n$ be
coordinates of point $p$ in some local chart in M, and let
$v^1,\,\ldots,\,v^n$ be components of vector $\bold v$
in this local chart:
$$
\hskip -2em
\bold v=v^1\,\frac{\partial}{\partial x^1}+\ldots
+v^n\,\frac{\partial}{\partial x^n}.
\tag1.1
$$
Lagrange function $L$ in local coordinates $x^1,\,\ldots,\,x^n$
is represented as the function of $2n$ arguments $L=L(x^1,,\ldots,
x^n,v^1,\ldots,v^n)$, while corresponding Lagrangian dynamical
system is given by the following ODE's:
$$
\xalignat 2
&\hskip -2em
\dot x^i=v^i,
&\frac{d}{dt}\!\left(\frac{\partial L}
{\partial v^i}\right)=\frac{\partial L}{\partial x^i}.
\tag1.2
\endxalignat
$$
Due to $\dot x^i=v^i$ in \thetag{1.2} vector \thetag{1.1} is called
{\it velocity vector}. In addition to \thetag{1.2} we consider so
called modified Lagrange equations:
$$
\xalignat 2
&\hskip -2em
\dot x^i=\frac{v^i}{\Omega},
&&\frac{d}{dt}\!\left(\frac{\partial L}{\partial v^i}\right)=
\frac{1}{\Omega}\,\frac{\partial L}{\partial x^i}.
\tag1.3
\endxalignat
$$
Denominator $\Omega$ in these equations \thetag{1.3} is determined by
the following formula:
$$
\hskip -2em
\Omega=\sum^n_{i=1}v^i\,\frac{\partial L}{\partial v^i}.
\tag1.4
$$
Like original Lagrange equations \thetag{1.2}, modified equations
\thetag{1.3} with denominator \thetag{1.4} arise in applications. As
shown in \cite{1}, they describe wave front dynamics for various wave
propagation phenomena. Now let's denote by $\bold p$ the covector at
the point $p\in M$ given by its components in local chart:
$$
\hskip -2em
p_i=\frac{\partial L}{\partial v^i}.
\tag1.5
$$
Covector $\bold p$ is called {\it momentum covector}. Pair composed
by point $p$ and covector $\bold p$ is a point of cotangent bundle
$T^*\!M$. Therefore \thetag{1.5} defines a map
$$
\hskip -2em
\lambda\!:TM\to T^*\!M.
\tag1.6
$$
This map is called {\it Legendre transformation}. It is well known in
mechanics (see \cite{2}).\par
    Note that first equations in \thetag{1.2} and \thetag{1.3}
are different. Therefore vector $\bold v$ cannot be interpreted
as velocity vector for modified dynamical system \thetag{1.3}.
For this purpose we introduce vector $\bold u$ with components
$$
\hskip -2em
u^i=\frac{v^i}{\Omega}.
\tag1.7
$$
Formula \thetag{1.7} determines another map, which is similar to
\thetag{1.6}:
$$
\hskip -2em
\mu\!:TM\to TM.
\tag1.8
$$
\definition{Definition 1.1} Lagrangian dynamical systems \thetag{1.2}
and \thetag{1.3} and their Lagrange function $L$ are called {\it regular}
if both maps \thetag{1.4} and \thetag{1.8} are diffeomorphisms and if
denominator $\Omega$ in modified Lagrange equations \thetag{1.3} given
by formula \thetag{1.4} is positive at all points $q=(p,\bold v)$ of
tangent bundle $TM$, where $|\bold v|\neq 0$.
\enddefinition
    In order to get simplest example of regular Lagrangian dynamical
system one should assume $M$ to be equipped with Riemannian metric
$\bold g$ and one should take Lagrange function $L$ to be quadratic
with respect to velocity vector $\bold v$:
$$
\hskip -2em
L=\frac{1}{2}\sum^n_{i=1}\sum^n_{j=1}g_{ij}\,v^i\,v^j-U(x^1,\ldots,x^n).
\tag1.9
$$
Though it is simple, this example covers all real mechanical systems
if we neglect friction in junctions and other forms of energy pumping
and dissipation. First term in \thetag{1.9} is kinetic energy, while
second term $U(x^1,\ldots,x^n)$ is potential energy of mechanical system.
For $L$ of the form \thetag{1.9} Legendre transformation \thetag{1.6}
looks like index lowering procedure in metric $\bold g$. Indeed, applying
\thetag{1.5} to \thetag{1.9}, we find
$$
\hskip -2em
p_i=\sum^n_{j=1}g_{ij}\,v^j.
\tag1.10
$$
If $L$ is treated as geometric equipment of manifold $M$, then for $L$
of the form \thetag{1.9} this equipment is equivalent to Riemannian metric
in essential. This is the case considered in papers \cite{3} and \cite{4}.
In present paper we consider more general case, when geometric equipment
of manifold $M$ is determined by some arbitrary regular Lagrange function
$L$, which is not necessarily given by formula \thetag{1.9}.\par
    Now let's consider another dynamical system. It can be either
Lagrangian or not Lagrangian, but we assume it to be second order
dynamical system. More precisely, we assume that it is given by
differential equations
$$
\xalignat 2
&\hskip -2em
\dot x^i=u^i,
&&\dot u^i=\Phi^i(x^1,\ldots,x^n,u^1,\ldots,u^n).
\tag1.11
\endxalignat
$$
We call \thetag{1.11} Newtonian dynamical system, since these equation
are similar to those expressing Newton's second law for the motion of
a particle of unit mass $m=1$ under the action of force $\boldsymbol
\Phi$. Certainly, we cannot touch all problems associated with pairs
of dynamical systems of the form \thetag{1.2} and \thetag{1.11}. In
present paper we construct theory of normal shift for such pairs of
dynamical systems, thus realizing the project claimed in previous
paper \cite{1}.\par
    Theory of Newtonian dynamical systems {\it admitting normal shift}
was initiated in 1993 in preprint \cite{5}. At first dynamical systems
\thetag{1.11} in Euclidean space $\Bbb R^n$ were considered, then theory
was extended for dynamical systems in Riemannian and Finslerian manifolds.
This phase, which lasted 7 years from 1993 till 1998, is reflected in
theses \cite{6} and \cite{7} (see also appropriate references therein).
For recent results in the theory of dynamical systems admitting normal
shift see papers \cite{8--12} and papers \cite{1}, \cite{3}, and
\cite{4} already mentioned above. Below in section~2 we give some
preliminary information and motivated definitions. Then in further
sections of this paper we construct theory of Newtonian dynamical
systems \thetag{1.11} admitting normal shift for the case of manifolds,
geometry of which is not Euclidean, not Riemannian, and even not
Finslerian, but is given by Lagrange function $L$ of some regular
Lagrangian dynamical system in them.
\head
2. Normal shift of hypersurfaces.
\endhead
\parshape 32 0pt 360pt 0pt 360pt
0pt 360pt 0pt 360pt 0pt 360pt 0pt 360pt 0pt 360pt
0pt 360pt 0pt 360pt 0pt 360pt 0pt 360pt 0pt 360pt 0pt 360pt
0pt 360pt 0pt 360pt 0pt 360pt 160pt 200pt 160pt 200pt 160pt 200pt
160pt 200pt 160pt 200pt 160pt 200pt 160pt 200pt 160pt 200pt 
160pt 200pt 160pt 200pt 160pt 200pt 160pt 200pt 160pt 200pt
160pt 200pt 160pt 200pt
0pt 360pt
     Suppose that $S$ is some arbitrary smooth hypersurface in $M$. We
say that $S$ is equipped with smooth transversal vector field if at each
point $p\in S$ some nonzero vector $\bold u(p)$ transversal to $S$ is
fixed. In this case one can consider the following initial data for
Newtonian dynamical system \thetag{1.11}:
$$
\xalignat 2
&\hskip -2em
x^i\,\hbox{\vrule height 8pt depth 8pt width 0.5pt}_{\,t=0}
=x^i(p),
&&u^i\,\hbox{\vrule height 8pt depth 8pt width 0.5pt}_{\,t=0}=
u^i(p).
\tag2.1
\endxalignat
$$
Here $x^i(p)$ are coordinates of point $p$ in some local chart in $M$
and $u^i(p)$ are components of transversal vector $\bold u(p)$ in this
chart. Applying initial data \thetag{2.1} to \thetag{1.11}, we get a
family of trajectories of this dynamical system starting at the points
of $S$. In local chart it is represented by functions
$$
\pagebreak
\hskip -2em
\cases x^1=x^1(t,p),\\ .\ .\ .\ .\ .\ .\ .\ .\ .\ \\
x^n=x^n(t,p).\endcases\hskip -1em
\tag2.2
$$
These trajectories are transversal to $S$. If we fix time instant $t\neq
0$ and gather all points of trajectories \thetag{2.2} corresponding to
this time instant (see Fig.~2.1), we get another hypersurface $S_t$ and
diffeomorphism
$$
f_t\!:S\to S_t
\tag2.3
$$
binding $S_t$ with initial hypersurface\footnote{In general, this is true
only for sufficiently small hypersurface $S$ and for time instants $t$
sufficiently close to initial time instant $t=0$. However, our further
consideration is local. Therefore here we shall not discuss the problem
of globalization for diffeomorphisms \thetag{2.3} referring reader to 
paper \cite{12}, where some aspects of this problem are studied}.
\vadjust{\vskip -36pt\hbox to 0pt{\kern 5pt\hbox{\special{em:graph
pst-14a.gif}}\hss}\vskip 36pt}Diffeo\-morphism $f_t$ (or, more precisely,
the whole set of diffeomorphisms
$f_t$) is called {\it a shift} of $S$ along trajectories of Newtonian
dynamical system \thetag{1.11}. Note that the shift \thetag{2.3} keeps
transversality in local. This means that trajectories of shift
are transversal not only to initial hypersurface $S$, but to all
hypersurfaces $S_t$ for sufficiently small values of $t$. Shift $f_t$ is
called {\it normal shift} if it keeps orthogonality of $S_t$ and
trajectories in some sense. In previous papers (see \cite{8--12} and
earlier) orthogonality was understood in the sense of some metric either
Euclidean, Riemannian, or Finslerian. Below we shall understand it in the
sense of Lagrange function $L$ as it was suggested in paper \cite{1}.\par
\adjustfootnotemark{-1}
    Let $y^1,\,\ldots,\,y^{n-1}$ be some local coordinates on initial
hypersurface $S$. Due to diffeomorphisms of shift \thetag{2.3} they can
be transferred to all hypersurfaces $S_t$. Functions \thetag{2.2} in terms
of local coordinates $y^1,\,\ldots,\,y^{n-1}$ are written as follows:
$$
\hskip -2em
\cases
x^1=x^1(t,y^1,\ldots,y^{n-1}),\\
.\ .\ .\ .\ .\ .\ .\ .\ .\ .\ .\ .\ .\ .\ .\ .\ .\ .\ .\\
x^n=x^n(t,y^1,\ldots,y^{n-1}).
\endcases
\tag2.4
$$
Their time derivatives $u^i=\dot x^i$ are components of velocity vector
$\bold u$. It's easy to understand that partial derivatives of these
functions with respect to $y^i$ are components of vector tangent to $S_t$.
It is called $i$-th vector of variation (or, more precisely, $i$-th vector
of variation of trajectories). We denote this vector by $\boldsymbol\tau_i$:
$$
\hskip -2em
\boldsymbol\tau_i=\tau^1_i\,\frac{\partial}{\partial x^1}+\ldots
+\tau^n_i\,\frac{\partial}{\partial x^n}\text{, \ where \ }
\tau^s_i=\frac{\partial x^s}{\partial y^i}
\tag2.5
$$
Vectors $\boldsymbol\tau_1,\,\ldots,\,\boldsymbol\tau_{n-1}$ form a base
in tangent space to $S_t$. For normal shift they should be perpendicular
to velocity vector $\bold u$. According to receipt from paper \cite{1},
we define the following deviation functions $\varphi_i$:
$$
\hskip -2em
\varphi_i=\left<\bold p\,|\,\boldsymbol\tau_i\right>=
\sum^n_{s=1}p_s\,\tau^s_i.
\tag2.6
$$
Here we have no metric, therefore scalar product $\left<\bold p\,|\,
\boldsymbol\tau_i\right>$ is nothing else, but symbolic notation for
the sum in right hand side of \thetag{2.6}.
\definition{Definition 2.1} Shift of hypersurface $S$ along trajectories
of Newtonian dynamical system \thetag{1.11} defined by initial data
\thetag{2.1} is called {\it normal shift} if all deviation functions
$\varphi_i$ given by formula \thetag{2.6} are identically zero.
\enddefinition
Deviation functions $\varphi_1,\,\ldots,\,\varphi_{n-1}$ are used
as a measure of deviation of shift $f_t$ from being a normal shift.
Their vanishing is indicator of normality.
\head
3. Relative form of the equations of Newtonian dynamics.
\endhead
    Note that covector $\bold p$ in \thetag{2.6}, according to
\thetag{1.5}, depend on components of vector $\bold v$. However,
$\bold v$ is not velocity vector for dynamical system \thetag{1.11}.
Vectors $\bold v$ and $\bold u$ are bound by the map \thetag{1.8},
which is expressed by formula \thetag{1.7} in local chart. This map
is diffeomorphism due to our assumptions (see definition~1.1).
Therefore we can transform \thetag{1.11} to variables
$x^1,\,\ldots,\,x^n,\,v^1,\,\ldots,\,v^n$. Here we write
$$
\xalignat 2
&\hskip -2em
\dot x^i=\frac{v^i}{\Omega},
&&\dot v^i=\Psi^i(x^1,\ldots,x^n,v^1,\ldots,v^n).
\tag3.1
\endxalignat
$$
One could express $\Psi^1,\,\ldots,\,\Psi^n$ through functions
$\Phi^1,\,\ldots,\,\Phi^n$ in \thetag{1.11}. However, the latter
ones are arbitrary functions. Therefore we can assume $\Psi^1,\,
\ldots,\,\Psi^n$ in \thetag{3.1} to be arbitrary functions as well,
with no need to follow their relations to the functions $\Phi^1,\,
\ldots,\,\Phi^n$ in \thetag{1.11}.\par
     In the next step we use Legendre transformation \thetag{1.6}
in order to transform differential equations \thetag{3.1} further.
Now we write them as
$$
\xalignat 2
&\hskip -2em
\dot x^i=\frac{v^i}{\Omega},
&&\frac{d}{dt}\!\left(\frac{\partial L}{\partial v^i}\right)
-\frac{1}{\Omega}\,\frac{\partial L}{\partial x^i}=Q_i,
\tag3.2
\endxalignat
$$
where $Q_i=Q_i(x^1,\ldots,x^n,v^1,\ldots,v^n)$. Since $\lambda$ is
diffeomorphism, these equation \thetag{3.2} are equivalent to
previous ones \thetag{3.1}. This is relative form of the equations
of Newtonian dynamics \thetag{1.11}. Relative, because in writing
them we need another dynamical system \thetag{1.3}.
\head
4. Extended tensor fields.
\endhead
    Let's consider quantities $Q_1,\,\ldots,\,Q_n$ in \thetag{3.2}.
By means of direct calculations one can check up that these
quantities are transformed as components of covector under the
changes of local charts in $M$. They define a covector at the point
$p$, where $p$ is a point with local coordinates $x^1,\,\ldots,\,x^n$.
Let's denote this covector by $\bold Q$ and note that it depends not
only on $p$, but on components of vector $\bold v$ as well. This means
that $\bold Q$ fits the definition of extended covector field.
\definition{Definition 4.1} Extended covector field $\bold X$ on
a manifold $M$ is a covector-valued function that to each point $q=(p,
\bold v)$ of tangent bundle $TM$ puts into correspondence some covector
from cotangent space $T^*_p(M)$ at the point $p=\pi(q)$ in $M$.
\enddefinition
    Similarly one can define extended tensor fields on the manifold
$M$. First let's consider the following tensor product of tangent
and cotangent spaces:
$$
T^r_s(p,M)=\overbrace{T_p(M)\otimes\ldots\otimes T_p(M)}^{\text{$r$
times}}\otimes\underbrace{T^*_p(M)\otimes\ldots\otimes T^*_p(M)}_{\text{$s$
times}}.
$$
Tensor product $T^r_s(p,M)$ is known as a space of $(r,s)$-tensors at the
point $p$ of the manifold $M$. Pair of integer numbers $(r,s)$ determines
the type of tensors. Elements of the space $T^r_s(p,M)$ are called
{\it $r$-times contravariant and $s$-times covariant tensors}, or tensors
of the type $(r,s)$, or, for brevity, $(r,s)$-tensors.
\definition{Definition 4.2} Extended tensor field $\bold X$
of the type $(r,s)$ on a manifold $M$ is a tensor-valued function
that to each point $q=(p,\bold v)$ of tangent bundle $TM$ puts into
correspondence some tensor from tensor space $T^r_s(p,M)$.
\enddefinition
    As far as I know, extended tensor fields first arose in Finslerian
geometry. They was intensively used in theory of Newtonian dynamical
systems admitting normal shift (see theses \cite{6}, \cite{7} and
references therein). Their application to Lagrangian and Hamiltonian
dynamical systems is explained in \cite{1}, \cite{3}, and in \cite{4}.
\head
5. Momentum representation for extended tensor fields.
\endhead
    Note that we can replace tangent bundle  $TM$ in definition~4.2
by cotangent bundle $T^*\!M$. Then we obtain another definition of
extended tensor field.
\definition{Definition 5.1} Extended tensor field $\bold Y$ of the
type $(r,s)$ on a manifold $M$ is a tensor-valued function that to
each point $q=(p,\bold p)$ of cotangent bundle $T^*\!M$ puts into
correspondence some tensor from tensor space $T^r_s(p,M)$.
\enddefinition
Extended tensor fields as given by definitions~4.2 and 5.1 are
two different objects. However, they can be related to each other
by Legendre map \thetag{1.6} that links tangent bundle $TM$ and
cotangent bundle $T^*\!M$. Suppose that
$$
\xalignat 2
&\hskip -2em
\bold X=\bold Y\compos\lambda,
&&\bold Y=\bold X\compos\lambda^{-1}.
\tag5.1
\endxalignat
$$
If relationships \thetag{5.1} hold, we say that $\bold Y$ is
$\bold p$-representation (or {\bf momentum} representation) of
extended tensor field $\bold X$. Similarly, in this case we
say that $\bold X$ is $\bold v$-representation (or {\bf velocity}
representation) of extended tensor field $\bold Y$.
\head
6. Weak normality condition.
\endhead
    Let's consider Newtonian dynamical system written in relative
form \thetag{3.2} and let's consider some one-parametric family of
trajectories $p=p(t,y)$ of this dynamical system. In local chart
this family of trajectories is represented by functions 
$$
\hskip -2em
\cases
x^1=x^1(t,y),\\
.\ .\ .\ .\ .\ .\ .\ .\ .\ .\\
x^n=x^n(t,y)
\endcases
\tag6.1
$$
(compare with \thetag{2.4}). Time derivatives of \thetag{6.1} define
vector $\bold v$ according to first part of the equations \thetag{3.2}.
Components of $\bold v$ depend on $t$ and $y$:
$$
\hskip -2em
v^1=v^1(t,y),\ \dots,\ v^n=v^n(t,y).
\tag6.2
$$
Like in \thetag{2.5}, we can define variation vector $\boldsymbol\tau$.
It is given by formula
$$
\hskip -2em
\boldsymbol\tau=\tau^1\,\frac{\partial}{\partial x^1}+\ldots
+\tau^n\,\frac{\partial}{\partial x^n}\text{, \ where \ }
\tau^s=\frac{\partial x^s}{\partial y}.
\tag6.3
$$
Here we have only one variation vector since, besides time variable
$t$, in \thetag{6.1} we have only one parameter $y$. Therefore we have
only one deviation function
$$
\hskip -2em
\varphi=\left<\bold p\,|\,\boldsymbol\tau\right>=\sum^n_{s=1}p_s
\,\tau^s.
\tag6.4
$$
Let's differentiate \thetag{6.2} with respect to parameter $y$.
This yields a series of functions
$$
\hskip -2em
\theta^i=\theta^i(t,y)=\frac{\partial v^i}{\partial y}.
\tag6.5
$$
Double set of functions $\tau^1,\,\ldots,\,\tau^n$ and $\theta^1,\,
\ldots,\,\theta^n$ satisfy a system of linear ordinary differential
equations with respect to time variable $t$. This system is obtained
as linearization for \thetag{3.2}. Differentiating first part of
the equations \thetag{3.2} with respect to parameter $y$, we get
the following expression for time derivative $\dot\tau^i$:
$$
\hskip -2em
\gathered
\dot\tau^i=\frac{\theta^i}{\Omega}-\sum^n_{s=1}\frac{v^i}
{\Omega^2}\,\frac{\partial L}{\partial v^s}\,\theta^s\,-\\
\vspace{1ex}
-\sum^n_{k=1}\sum^n_{s=1}\frac{\partial^2\!L}{\partial v^k\,
\partial v^s}\,\frac{v^i\,v^k\,\theta^s}{\Omega^2}
-\sum^n_{k=1}\sum^n_{s=1}\frac{\partial^2\!L}{\partial v^k\,
\partial x^s}\,\frac{v^i\,v^k\,\tau^s}{\Omega^2}.
\endgathered
\tag6.6
$$
Differentiating second part of the equations \thetag{3.2}
with respect to $y$, we obtain
$$
\hskip -2em
\gathered
\sum^n_{s=1}\frac{\partial^2\!L}{\partial v^i\,\partial v^s}\,
\dot\theta^s+\sum^n_{s=1}\frac{\partial^2\!L}{\partial v^i\,
\partial x^s}\,\dot\tau^s+
\sum^n_{k=1}\sum^n_{s=1}\frac{\partial^3\!L}{\partial v^i\,
\partial v^s\,\partial v^k}\,\dot v^k\,\theta^s\,+\\
+\sum^n_{k=1}\sum^n_{s=1}\frac{\partial^3\!L}{\partial v^i\,
\partial v^s\,\partial x^k}\,\dot x^k\,\theta^s
+\sum^n_{k=1}\sum^n_{s=1}\frac{\partial^3\!L}{\partial v^i\,
\partial x^s\,\partial v^k}\,\dot v^k\,\tau^s\,+\\
\vspace{1ex}
+\sum^n_{k=1}\sum^n_{s=1}\frac{\partial^3\!L}{\partial v^i\,
\partial x^s\,\partial x^k}\,\dot x^k\,\tau^s
-\sum^n_{s=1}\frac{\partial^2\!L}{\partial x^i\,
\partial v^s}\,\frac{\theta^s}{\Omega}-\sum^n_{s=1}
\frac{\partial^2\!L}{\partial x^i\,\partial x^s}\,\frac{\tau^s}{\Omega}\,+\\
\vspace{1ex}
+\sum^n_{s=1}\frac{\partial L}{\partial x^i}\,
\frac{\partial L}{\partial v^s}\,\frac{\theta^s}{\Omega^2}
+\sum^n_{k=1}\sum^n_{s=1}\frac{\partial L}{\partial x^i}\,
\frac{\partial^2\!L}{\partial v^k\,\partial v^s}
\frac{\,v^k\,\theta^s}{\Omega^2}\,+\\
\vspace{1ex}
+\sum^n_{k=1}\sum^n_{s=1}\frac{\partial L}{\partial x^i}\,
\frac{\partial^2\!L}{\partial v^k\,\partial x^s}
\frac{\,v^k\,\tau^s}{\Omega^2}=\sum^n_{s=1}\frac{\partial Q_i}
{\partial x^s}\,\tau^s+\sum^n_{s=1}\frac{\partial Q_i}
{\partial v^s}\,\theta^s.
\endgathered
\tag6.7
$$
The equalities \thetag{6.7} are not resolved with respect to
derivatives $\dot\theta^1,\,\ldots,\,\dot\theta^n$. However, they
can be resolved. Indeed, if we denote by $\bold g$ a matrix with
components
$$
\hskip -2em
g_{ij}=\frac{\partial^2\!L}{\partial v^i\,\partial v^j},
\tag6.8
$$
we can see that due to \thetag{1.5} it is Jacoby matrix for Legendre
transformation \thetag{1.6}. Matrix \thetag{6.8} is non-degenerate
since Legendre map $\lambda$ is diffeomorphism due to our initial
assumptions (see definition~1.1 and ending part of section 1 after
it).\par
    As we already mentioned above, both \thetag{6.6} and \thetag{6.7}
form a system of homogeneous linear ordinary differential equations
with respect to functions $\tau^1,\,\ldots,\,\tau^n$ and $\theta^1,\,
\ldots,\,\theta^n$. Solutions of these equations form linear space
of dimension $2n$, We denote it by $\goth T$. Looking at \thetag{6.4},
one can see that $\varphi$ depends linearly on $\tau^1,\,\ldots,\,
\tau^n$. Hence all time derivatives of $\varphi$ (i\.\,e\. $\dot\varphi$,
$\ddot\varphi$, and so forth) depend linearly on $\tau^1,\,\ldots,\,
\tau^n$ and on time derivatives of $\tau^1,\,\ldots,\,\tau^n$. Due to
equations \thetag{6.6} and \thetag{6.7} the latter ones can be expressed
linearly through $\tau^1,\,\ldots,\,\tau^n$ and $\theta^1,\,\ldots,
\,\theta^n$. This means that for each fixed instance of time $t$ the
value of function $\varphi(t)$ itself and values of all time derivatives
of this function
$$
\varphi^{(k)}=\frac{d^k\!\varphi}{dt^k}
$$
are linear functionals belonging to dual space $\goth T^*$. The dimension
of $\goth T^*$ is finite, it is equal to $2n$. Therefore functions
$\varphi,\ \dot\varphi,\ \ddot\varphi,\ \ldots,\ \varphi^{(2n)}$
are linearly dependent as elements of $\goth T^*$. This means that
there are some coefficients $C_0,\,\ldots,\,C_{2n}$, which do not vanish
simultaneously, such that the following equality holds:
$$
\hskip -2em
\sum^n_{i=0}C_i\,\varphi^{(i)}=0.
\tag6.9
$$
For fixed $t$ coefficients $C_0,\,\ldots,\,C_{2n}$ in \thetag{6.9} are
real numbers. However, if $t$ is not fixed, then $C_0,\,\ldots,\,C_{2n}$
depend on $t$. They also depend on that particular trajectory of dynamical
system \thetag{3.2}, for which functions $\varphi^{(i)}(t)$ are
calculated:
$$
\hskip -2em
\sum^n_{i=0}C_i(t)\,\varphi^{(i)}=0.
\tag6.10
$$
\proclaim{Theorem 6.1} For each trajectory $p=p(t)$ of Newtonian
dynamical system \thetag{3.2} all deviation functions $\varphi=
\varphi(t)$ on this trajectory satisfy the same linear homogeneous
ordinary differential equation \thetag{6.10} of the order not
greater than $2n$.
\endproclaim
   Saying {\it ``all deviation functions''}, in theorem~6.1 we assume
that given trajectory $p=p(t)$ can be included into one-parametric
family of trajectories by all possible ways. Each such inclusion
defines some variation vector $\boldsymbol\tau=\boldsymbol\tau(t)$
and corresponding deviation function $\varphi=\varphi(t)$ on that
trajectory. Functions $C_0,\,\ldots,\,C_{2n}$ in \thetag{6.10}
depend on the trajectory $p=p(t)$, but they do not depend on how
this trajectory is included into one parametric family of
trajectories.\par
   Theorem~6.1 gives upper estimate for the order of ODE \thetag{6.10}.
For most cases this estimate $2n$ is reached. However, in some special
cases real order can be much less than $2n$. We consider one of such
special cases.
\definition{Definition 6.1} We say that Newtonian dynamical system
\thetag{3.2} satisfies {\bf weak} normality condition if for each
its trajectory $p=p(t)$ there is some second order homogeneous
linear ordinary differential equation 
$$
\hskip -2em
\ddot\varphi=\Cal A(t)\,\dot\varphi+\Cal B(t)\,\varphi
\tag6.11
$$
such that all deviation functions on this trajectory satisfy this
differential \pagebreak equation.
\enddefinition
\noindent In other words, weak normality condition is fulfilled
when \thetag{6.10} reduces to second order differential equation
\thetag{6.11}.
\head
7. Additional normality condition.
\endhead
    Suppose that weak normality condition is fulfilled somehow. Then
all deviation functions \thetag{2.6} arising in shift construction
satisfy second order differential equation of the form \thetag{6.11}.
Coefficients $\Cal A$ and $\Cal B$ depend on shift trajectory (see
definition~6.1), therefore here we should write this equation as
$$
\hskip -2em
\ddot\varphi=\Cal A(t,y^1,\ldots,y^{n-1})\,\dot\varphi
+\Cal B(t,y^1,\ldots,y^{n-1})\,\varphi.
\tag7.1
$$
In spite of this minor difference \thetag{7.1} is linear homogeneous
ODE for $\varphi$ with respect to time variable $t$. According to
definition~2.1, in order to have normal shift we should provide
vanishing of all deviation functions $\varphi_i,\,\ldots,\,\varphi_{n-1}$
in \thetag{2.6}. Due to \thetag{7.1} now it is sufficient to provide
the following initial data for them:
$$
\xalignat 2
&\hskip -2em
\varphi_i,\hbox{\vrule height 8pt depth 8pt width 0.5pt}_{\,t=0}=0,
&&\dot\varphi_i\,\hbox{\vrule height 8pt depth 8pt width 0.5pt}_{\,t=0}=0.
\tag7.2
\endxalignat
$$\par
The only way to provide these initial conditions is to choose proper
initial conditions in \thetag{2.1}. First part of initial conditions
\thetag{2.1} says that shift trajectories should start the points of
initial hypersurface $S$. We cannot change this condition. However, we
can specify second part of initial conditions \thetag{2.1}. Note that
at each point $p\in S$ we have tangent space $T_p(S)$ embedded into
tangent space $T_p(M)$ as a subspace of codimension $1$. Denote by
$\bold n=\bold n(p)$ some nonzero covector that vanishes when contracted
with all vectors from $T_p(S)$:
$$
\hskip -2em
\left<\bold n\,|\,\boldsymbol\tau\right>=\sum^n_{s=1}n_s\,\tau^s=0
\text{, \ for all \ }\boldsymbol\tau\in T_p(S).
\tag7.3
$$
It is called {\it normal covector} for $S$. At each point $p\in S$ normal
covector $\bold n$ is determined uniquely up to a scalar factor. We have
no tools for to specify this factor canonically, but, nevertheless, we can
glue normal covectors from various points into a smooth covector-valued
function $\bold n=\bold n(p)$ on $S$ (or locally in small domains covering
all points of $S$). Now let's compare \thetag{7.3} with \thetag{2.6}
and remember that for $t=0$ variation vectors $\boldsymbol\tau_1,\,\ldots,
\,\boldsymbol\tau_{n-1}$ form a base in tangent spaces to $S$. Therefore
in order to provide first part of initial conditions \thetag{7.1} we
should direct initial momentum covector $\bold p$ along normal covector
of $S$. This means that we should specify initial conditions \thetag{2.1}
as follows:
$$
\xalignat 2
&\hskip -2em
x^i\,\hbox{\vrule height 8pt depth 8pt width 0.5pt}_{\,t=0}
=x^i(p),
&&p_i\,\hbox{\vrule height 8pt depth 8pt width 0.5pt}_{\,t=0}=
\nu(p)\cdot n_i(p).
\tag7.4
\endxalignat
$$
Here $n_i(p)$ are components of normal covector $\bold n(p)$, while
$\nu=\nu(p)$ is some smooth scalar function on $S$, which is yet
undefined.\par
   When applied to Newtonian dynamical system written in relative
form \thetag{3.2}, initial data \thetag{7.4} determine initial
velocity $\bold v$ implicitly through Legendre transformation
\thetag{1.5}. \pagebreak However, we can make them explicit initial data if
we pass to $\bold p$-representation in the equations \thetag{3.2}.
Here these equations look like
$$
\xalignat 2
&\hskip -2em
\dot x^i=\frac{1}{\Omega}\,\frac{\partial H}{\partial p_i},
&&\dot p_i=-\frac{1}{\Omega}\,\frac{\partial H}{\partial x^i}
+Q_i.
\tag7.5
\endxalignat
$$
Equations \thetag{7.5} are the equations of Newtonian dynamical
system in the form relative to modified Hamiltonian dynamical
system with Hamilton function $H$. Function $H$ is derived from
Lagrange function $L$ as follows:
$$
\hskip -2em
H=h\compos\lambda^{-1}\text{, \ where \ }
h=\sum^n_{i=1}v^i\,\frac{\partial L}{\partial v^i}-L.
\tag7.6
$$
In other words, extended scalar field $H$ is $\bold p$-representation
of extended scalar field $h$, where $h$ is given by formula
\thetag{7.6}. Denominator $\Omega$ in \thetag{7.5} is
$\bold p$-representation of denominator $\Omega$ in \thetag{3.2}.
It is expressed through Hamilton function $H$:
$$
\hskip -2em
\Omega=\sum^n_{i=1}p_i\,\frac{\partial H}{\partial p_i}.
\tag7.7
$$
Functions $Q_i$ in \thetag{7.5} are components of the same extended
covector field $\bold Q$, as in \thetag{3.2}, but transformed to
$\bold p$-representation. Their arguments are $x^1,\,\ldots,\,x^n$
and $p_1,\,\ldots,\,p_n$. More details concerning formulas
\thetag{7.6} and \thetag{7.7} and Legendre transformation in whole
can be found in book \cite{2} and in papers \cite{1}, \cite{3}, and
\cite{4}.\par
    Now let's proceed with initial conditions \thetag{7.2}. First
part of these initial conditions now is fulfilled due to \thetag{7.4}.
We should provide second part of them by proper choice of function
$\nu=\nu(p)=\nu(y^1,\ldots,y^{n-1})$ in \thetag{7.4}. Let's calculate
time derivative of deviation function $\varphi_i$ using differential
equations \thetag{7.5}:
$$
\hskip -2em
\gathered
\dot\varphi_i=\frac{d}{dt}\!\left(\,\shave{\sum^n_{s=1}}p_s\,
\tau^s_i\!\right)=\sum^n_{s=1}\dot p_s\,\tau^s_i
+\sum^n_{s=1}p_s\,\dot\tau^s_i=\\
\vspace{1ex}
=-\sum^n_{s=1}\frac{1}{\Omega}\,\frac{\partial H}{\partial x^s}\,
\tau^s_i+\sum^n_{s=1}Q_s\,\tau^s_i+\sum^n_{s=1}p_s\,\dot\tau^s_i.
\endgathered
\tag7.8
$$
In order to calculate time derivatives $\dot\tau^s_i$ in formula
\thetag{7.8} we also use differential equations \thetag{7.5}. As
a result we get the following formula for $\dot\tau^s_i$:
$$
\gather
\dot\tau^s_i=\frac{\partial^2 x^s}{\partial t\,\partial y^i}=
\frac{\partial}{\partial y^i}\!\left(\frac{1}{\Omega}\,\frac{\partial H}
{\partial p_s}\right)=
\sum^n_{r=1}\frac{\partial x^r}{\partial y^i}\,
\frac{\partial}{\partial x^r}\!\left(\frac{1}{\Omega}\,\frac{\partial H}
{\partial p_s}\right)\,+\\
\vspace{1ex}
+\sum^n_{r=1}\frac{\partial p_r}{\partial y^i}\,
\frac{\partial}{\partial p_r}\!\left(\frac{1}{\Omega}\,\frac{\partial H}
{\partial p_s}\right)=
\sum^n_{r=1}\tau^r_i\,\frac{\partial}{\partial x^r}\!\left(\frac{1}
{\Omega}\,\frac{\partial H}{\partial p_s}\right)
+\sum^n_{r=1}\frac{\partial p_r}{\partial y^i}\,
\frac{\partial}{\partial p_r}\!\left(\frac{1}{\Omega}\,\frac{\partial H}
{\partial p_s}\right).
\endgather
$$
Let's substitute this formula into \thetag{7.8}. Then for time
derivative $\dot\varphi_i$ we obtain
$$
\gather
\dot\varphi_i=-\sum^n_{s=1}\frac{1}{\Omega}\,
\frac{\partial H}{\partial x^s}\,\tau^s_i
+\sum^n_{s=1}Q_s\,\tau^s_i+\sum^n_{s=1}\sum^n_{r=1}
p_s\,\tau^r_i\,\frac{\partial}{\partial x^r}\!\left(\frac{1}
{\Omega}\,\frac{\partial H}{\partial p_s}\right)\,+\\
\displaybreak
+\sum^n_{s=1}\sum^n_{r=1}p_s\,\frac{\partial p_r}{\partial y^i}\,
\frac{\partial}{\partial p_r}\!\left(\frac{1}{\Omega}\,\frac{\partial H}
{\partial p_s}\right)=\sum^n_{s=1}Q_s\,\tau^s_i+\sum^n_{r=1}
\tau^r_i\,\frac{\partial}{\partial x^r}\!\left(\,\shave{\sum^n_{s=1}}
\frac{p_s}{\Omega}\,\frac{\partial H}{\partial p_s}\right)\,+\\
\vspace{1ex}
+\sum^n_{r=1}\frac{\partial p_r}{\partial y^i}\,
\frac{\partial}{\partial p_r}\!\left(\,\shave{\sum^n_{s=1}}\frac{p_s}
{\Omega}\,\frac{\partial H}{\partial p_s}\right)
-\sum^n_{s=1}\frac{1}{\Omega}\,\frac{\partial H}{\partial p_s}\,
\frac{\partial p_s}{\partial y^i}-\sum^n_{s=1}\frac{1}{\Omega}\,
\frac{\partial H}{\partial x^s}\,\tau^s_i.
\endgather
$$
Note that sums in round brackets are identically equal to unity:
$$
\hskip -2em
\sum^n_{s=1}\frac{p_s}{\Omega}\,\frac{\partial H}{\partial p_s}=1.
\tag7.9
$$
This follows from \thetag{7.7}. Due to \thetag{7.9} two terms in the
above formula for $\dot\varphi_i$ do vanish. And we get rather simple
formula for time derivative $\dot\varphi_i$:
$$
\hskip -2em
\dot\varphi_i=-\sum^n_{s=1}\frac{1}{\Omega}\,\frac{\partial H}
{\partial p_s}\,\frac{\partial p_s}{\partial y^i}
-\sum^n_{s=1}\frac{1}{\Omega}\,\frac{\partial H}{\partial x^s}
\,\tau^s_i+\sum^n_{s=1}Q_s\,\tau^s_i.
\tag7.10
$$
If we recall initial conditions \thetag{7.2}, then from \thetag{7.10}
we derive
$$
\hskip -2em
\sum^n_{s=1}\frac{1}{\Omega}\,\frac{\partial H}{\partial x^s}
\,\tau^s_i+\sum^n_{s=1}\frac{1}{\Omega}\,\frac{\partial H}
{\partial p_s}\,\left(\frac{\partial p_s}{\partial y^i}\right)\!
\,\hbox{\vrule height 12pt depth 8pt width 0.5pt}_{\,t=0}
=\sum^n_{s=1}Q_s\,\tau^s_i.
\tag7.11
$$
Calculating partial derivatives $\partial p_s/\partial y^i$ in
the equality \thetag{7.11}, we should remember that $p_s=\nu\cdot
n_s$. This follows from initial data \thetag{7.4}. Then
$$
\hskip -2em
\left(\frac{\partial p_s}{\partial y^i}\right)\!
\,\hbox{\vrule height 12pt depth 8pt width 0.5pt}_{\,t=0}
=\frac{\partial\nu}{\partial y^i}\,n_s+
\nu\,\frac{\partial n_s}{\partial y^i}
=\frac{1}{\nu}\,\frac{\partial\nu}{\partial y^i}\,p_s+
\nu\,\frac{\partial n_s}{\partial y^i}.
\tag7.12
$$
Substituting this expression into \thetag{7.11} and using formula
\thetag{7.7} for $\Omega$, we can transform \thetag{7.12} into the
partial differential equations for $\nu$:
$$
\hskip -2em
\frac{\partial\nu}{\partial y^i}=
-\sum^n_{s=1}\frac{\nu^2}{\Omega}\,\frac{\partial n_s}
{\partial y^i}\,\frac{\partial H}{\partial p_s}
-\sum^n_{s=1}\nu\left(\frac{1}{\Omega}\frac{\partial H}
{\partial x^s}-Q_s\!\right)\tau^s_i.
\tag7.13
$$
Having derived \thetag{7.13}, we have proved the following lemma.
\proclaim{Lemma 7.1}Initial conditions \thetag{7.2} for deviation
functions $\varphi_1,\,\ldots,\,\varphi_{n-1}$ are equivalent to
initial data \thetag{7.4} for shift trajectories, where scalar
function $\nu$ on $S$ satisfies differential equations \thetag{7.13}.
\endproclaim
    In two-dimensional manifold $M$, when $n=2$, hypersurfaces are curves.
In this case we have only one deviation function and only one variable
$y=y^1$ as a parameter on $S$. Therefore \thetag{7.13} appears to be
ordinary differential equation for the function $\nu=\nu(y)$. We can set
initial value problem
$$
\hskip -2em
\nu(y)\,\hbox{\vrule height 8pt depth 8pt width 0.5pt}_{\,y=0}=\nu_0,
\tag7.14
$$
which is always solvable (at least locally) for all $\nu_0\neq 0$.
\pagebreak This
means that in two-dimensional case weak normality condition stated in
definition~6.1 is sufficient for to arrange normal shift of any predefined
hypersurface $S$ in $M$.\par
    In multidimensional case $n\geqslant 3$ situation changes crucially.
Here equations \thetag{7.13} form so called {\it complete system of Pfaff
equation} for $\nu=\nu(y^1,\ldots,y^{n-1})$. Each separate equation in
such system can be treated as ODE. Therefore initial condition like
\thetag{7.14} is the best way for fixing some particular solution:
$$
\hskip -2em
\nu(p_0)=\nu_0.
\tag7.15
$$
However, initial value problem \thetag{7.15} for Pfaff equations
\thetag{7.13} is not always solvable: some additional conditions
should be fulfilled. This is why we had a fork in development of
theory of metric normal shift in Riemannian geometry (compare theses
\cite{6} and \cite{7}). This fork is present here as well.
\definition{Definition 7.1} Complete system of Pfaff equations
\thetag{7.13} is called {\it compatible} if initial value problem
\thetag{7.15} for it is locally solvable for all $\nu_0\neq 0$.
\enddefinition
\noindent Let's write Pfaff equations \thetag{7.13} formally,
denoting by $\psi_i$ their right hand sides:
$$
\hskip -2em
\frac{\partial\nu}{\partial y^i}=\psi_i(\nu,y^1,\ldots,y^{n-1})
\tag7.16
$$
Due to \thetag{7.16} we can calculate mixed partial derivatives
of $\nu$ in two different ways
$$
\align
\hskip -2em\frac{\partial^2\nu}{\partial y^i\,\partial y^j}&=
\frac{\partial\psi_i}{\partial y^j}+\frac{\partial\psi_i}{\partial\nu}
\,\psi_j=\vartheta_{ij}(\nu,y^1,\ldots,y^{n-1}),\hskip -2em
\tag7.17\\
\vspace{1ex}
\hskip -2em\frac{\partial^2\nu}{\partial y^j\,\partial y^i}&=
\frac{\partial\psi_j}{\partial y^i}+\frac{\partial\psi_j}{\partial\nu}
\,\psi_i=\vartheta_{ji}(\nu,y^1,\ldots,y^{n-1}).\hskip -2em
\tag7.18
\endalign
$$
Equating \thetag{7.17} and \thetag{7.18}, we get compatibility condition
for \thetag{7.16}:
$$
\hskip -2em
\vartheta_{ij}(\nu,y^1,\ldots,y^{n-1})=\vartheta_{ji}(\nu,y^1,
\ldots,y^{n-1}).
\tag7.19
$$
\proclaim{Lemma 7.2} Pfaff equations \thetag{7.16} are compatible in
the sense of definition~7.1 if and only if for $\nu\neq 0$ left and
right hands sides of \thetag{7.19} are equal to each other identically
as functions of $n$ independent variables $y^1,\,\ldots,\,y^{n-1}$,
and $\nu$.
\endproclaim
Lemma~7.2 is standard result in the theory of Pfaff equations. Proof
of this lemma can be found in thesis \cite{6}.
\definition{Definition 7.2} We say that Newtonian dynamical system
\thetag{3.2} or, equivalently, dynamical system \thetag{7.5} satisfies
{\bf additional} normality condition if Pfaff equations \thetag{7.13}
derived from initial conditions \thetag{7.2} are compatible for any
hypersurface $S$ in $M$ and for any marked point $p_0\in S$.
\enddefinition
\head
8. Complete and strong normality conditions.
\endhead
    Both {\bf weak} and {\bf additional} normality conditions constitute
so called {\bf complete} normality condition in Riemannian geometry (see
thesis \cite{6}). We shall keep this terminology saying that Newtonian
dynamical system satisfies {\bf complete} \pagebreak normality condition
if both conditions from definitions~6.2 and 7.2 are
fulfilled\footnote{In two-dimensional case $n=2$ complete normality
condition reduces to weak normality condition since additional normality
condition in this case is always fulfilled.}.
But, apart from this complete normality condition, we shall consider
{\bf strong} normality condition given by the following definition.
\adjustfootnotemark{-1}
\definition{Definition 8.1} Newtonian dynamical system given by
differential equations \thetag{3.2} or, equivalently, by differential
equations \thetag{7.5} satisfies {\bf strong} normality condition
if for any hypersurface $S$ in $M$, for any marked point $p_0\in S$,
and for any real constant $\nu_0\neq 0$ there is some smaller part
$S'$ of $S$ containing marked point $p_0$ and there is some smooth
function $\nu=\nu(p)$ in this smaller part $S'$ normalized by the
condition \thetag{7.15} and such that initial data \thetag{7.4} with
this function $\nu$ determine normal shift of $S'$ in the sense of
definition~2.1.
\enddefinition
   In simpler words, dynamical systems satisfying strong normality
condition are called systems {\bf admitting normal shift} of
hypersurfaces. They are able to move normally any predefined
hypersurface $S$ in $M$.\par
   It is easy to note that complete normality condition is sufficient
for Newtonian dynamical system to satisfy strong normality condition.
Below in section~17 we shall prove that it is not only sufficient, but
necessary condition as well.
\head
9. Compatibility equations.
\endhead
    Now suppose that $M$ is a manifold of dimension $n\geqslant 3$.
In this case we should study compatibility equations \thetag{7.19}
for Pfaff system \thetag{7.13}. For this purpose let's calculate
partial derivatives \thetag{7.17} and \thetag{7.18} explicitly:
$$
\gather
\frac{\partial^2\nu}{\partial y^i\,\partial y^j}=-\sum^n_{s=1}
\frac{\nu^2}{\Omega}\,\frac{\partial H}{\partial p_s}\,
\frac{\partial^2 n_s}{\partial y^i\,\partial y^j}-\sum^n_{s=1}
\nu\left(\frac{1}{\Omega}\,\frac{\partial H}{\partial x^s}-Q_s
\!\right)\frac{\partial^2 x^s}{\partial y^i\,\partial y^j}\,+
\vspace{1ex}
+\sum^n_{s=1}\sum^n_{r=1}\left(\frac{2\,\nu^3}{\Omega^2}\,
\frac{\partial H}{\partial p_s}
\,\frac{\partial H}{\partial p_r}\,\frac{\partial n_s}
{\partial y^i}\,\frac{\partial n_r}{\partial y^j}
+\nu\left(\frac{1}{\Omega}\,\frac{\partial H}
{\partial x^s}-Q_s\!\right)\!\!\left(\frac{1}{\Omega}\,
\frac{\partial H}{\partial x^r}-Q_r\!\right)\tau^s_i\,\tau^r_j\,+
\right.\\
\vspace{1ex}
\left.+\,\frac{2\,\nu^2}{\Omega}\left(\frac{1}{\Omega}\,
\frac{\partial H}{\partial x^s}-Q_s\!\right)\frac{\partial H}
{\partial p_r}\,\frac{\partial n_r}{\partial y^j}\,\tau^s_i
+\frac{\nu^2}{\Omega}\frac{\partial H}{\partial p_s}
\left(\frac{1}{\Omega}\,\frac{\partial H}{\partial x^r}-Q_r\!
\right)\frac{\partial n_s}{\partial y^i}\,\tau^r_j\right)-\\
\vspace{1ex}
-\sum^n_{s=1}\sum^n_{r=1}\nu^2\,\frac{\partial}{\partial p_s}\!
\left(\frac{1}{\Omega}\,\frac{\partial H}{\partial p_r}\right)
\frac{\partial p_s}{\partial y^i}\,\frac{\partial n_r}{\partial y^j}
-\sum^n_{s=1}\sum^n_{r=1}\nu\,\,\frac{\partial}{\partial p_s}\!
\left(\frac{1}{\Omega}\,\frac{\partial H}{\partial x^r}
-Q_r\!\right)
\frac{\partial p_s}{\partial y^i}\,\tau^r_j\,-\\
\vspace{1ex}
-\sum^n_{s=1}\sum^n_{r=1}\nu^2\,\frac{\partial}{\partial x^s}\!
\left(\frac{1}{\Omega}\,\frac{\partial H}{\partial p_r}\right)
\frac{\partial n_r}{\partial y^j}\,\tau^s_i
-\sum^n_{s=1}\sum^n_{r=1}\nu\,\,\frac{\partial}{\partial x^s}\!
\left(\frac{1}{\Omega}\,\frac{\partial H}{\partial x^r}
-Q_r\!\right)\tau^s_i\,\tau^r_j.
\endgather
$$
For Pfaff equations \thetag{7.13} to be compatible, right hand side
of the above equality should be symmetric in indices $i$ and $j$.
Some terms there are obviously symmetric. Below we shall not write
such terms explicitly denoting them by dots:
$$
\gather
\frac{\partial^2\nu}{\partial y^i\,\partial y^j}=\dots
+\sum^n_{s=1}\sum^n_{r=1}\left(\frac{2\,\nu^2}{\Omega}
\left(\frac{1}{\Omega}\,\frac{\partial H}{\partial x^s}
-Q_s\!\right)\frac{\partial H}{\partial p_r}
-\nu^2\,\frac{\partial}{\partial x^s}\!\left(\frac{1}{\Omega}
\,\frac{\partial H}{\partial p_r}\right)\right)
\frac{\partial n_r}{\partial y^j}\,\tau^s_i\,+\\
\displaybreak
+\sum^n_{s=1}\sum^n_{r=1}\left(\frac{\nu^2}{\Omega}\,
\frac{\partial H}{\partial p_s}\left(\frac{1}{\Omega}\,
\frac{\partial H}{\partial x^r}-Q_r\!\right)
-\nu^2\,\frac{\partial}{\partial p_s}\!\left(\frac{1}{\Omega}
\,\frac{\partial H}{\partial x^r}-Q_r\!\right)\right)
\frac{\partial n_s}{\partial y^i}\,\tau^r_j\,-\\
\vspace{1ex}
-\sum^n_{s=1}\sum^n_{r=1}\nu\,p_s\,\frac{\partial}{\partial p_s}\!
\left(\frac{1}{\Omega}\,\frac{\partial H}{\partial p_r}\right)
\frac{\partial\nu}{\partial y^i}\,\frac{\partial n_r}{\partial y^j}
-\sum^n_{s=1}\sum^n_{r=1}p_s\,\frac{\partial}{\partial p_s}\!
\left(\frac{1}{\Omega}\,\frac{\partial H}{\partial x^r}
-Q_r\!\right)\frac{\partial\nu}{\partial y^i}\,\tau^r_j\,-\\
\vspace{1ex}
-\sum^n_{s=1}\sum^n_{r=1}\nu^3\,\frac{\partial}{\partial p_s}\!
\left(\frac{1}{\Omega}\,\frac{\partial H}{\partial p_r}\right)
\frac{\partial n_s}{\partial y^i}\,\frac{\partial n_r}{\partial y^j}
-\sum^n_{s=1}\sum^n_{r=1}
\nu\,\frac{\partial}{\partial x^s}\!
\left(\frac{1}{\Omega}\,\frac{\partial H}{\partial x^r}
-Q_r\!\right)\tau^s_i\,\tau^r_j.
\endgather
$$
In the above calculations we used formula \thetag{7.12} in order to
express partial derivatives $\partial p_s/\partial y^i$ through
$\partial n_s/\partial y^i$. As a result we have got partial
derivatives $\partial\nu/\partial y^i$ in the above expression. In
order to eliminate them now we shall use \thetag{7.13}:
$$
\gather
\frac{\partial^2\nu}{\partial y^i\,\partial y^j}=\dots
+\sum^n_{s=1}\sum^n_{r=1}\left(\frac{2\,\nu^2}{\Omega}
\left(\frac{1}{\Omega}\,\frac{\partial H}{\partial x^s}
-Q_s\!\right)\frac{\partial H}{\partial p_r}
-\nu^2\,\frac{\partial}{\partial x^s}\!\left(\frac{1}
{\Omega}\,\frac{\partial H}{\partial p_r}\right)\,+
\vphantom{\shave{\sum^n_{q=1}}}\right.\\
\vspace{1ex}
\left.+\shave{\sum^n_{q=1}}\nu^2\,p_q\,\frac{\partial}
{\partial p_q}\!\left(\frac{1}{\Omega}\,\frac{\partial H}
{\partial p_r}\right)\!\left(\frac{1}{\Omega}\,\frac{\partial H}
{\partial x^s}-Q_s\!\right)\right)
\frac{\partial n_r}{\partial y^j}\,\tau^s_i
+\sum^n_{s=1}\sum^n_{r=1}\left(\frac{\nu^2}{\Omega}\,
\frac{\partial H}{\partial p_s}\,\times
\vphantom{\shave{\sum^n_{q=1}}}\right.\\
\vspace{1ex}
\left.\times
\left(\frac{1}{\Omega}\,\frac{\partial H}{\partial x^r}
-Q_r\!\right)-\nu^2\,\frac{\partial}{\partial p_s}\!
\left(\frac{1}{\Omega}\,\frac{\partial H}{\partial x^r}
-Q_r\!\right)
+\shave{\sum^n_{q=1}}\frac{\nu^2\,p_q}{\Omega}\,
\frac{\partial}{\partial p_q}\!\left(\frac{1}{\Omega}\,
\frac{\partial H}{\partial x^r}-Q_r\!\right)
\times\right.\\
\vspace{1ex}
\left.\vphantom{\shave{\sum^n_{q=1}}}
\times\frac{\partial H}{\partial p_s}\right)
\frac{\partial n_s}{\partial y^i}\,\tau^r_j
+\sum^n_{s=1}\sum^n_{r=1}\left(\,\shave{\sum^n_{q=1}}
\frac{\nu^3\,p_q}{\Omega}\,\frac{\partial}{\partial p_q}\!
\left(\frac{1}{\Omega}\,\frac{\partial H}{\partial p_r}\right)
\frac{\partial H}{\partial p_s}+\frac{\nu^3}{\Omega^2}\,  
\frac{\partial\Omega}{\partial p_s}\,\frac{\partial H}
{\partial p_r}\right)\times\\
\vspace{1ex}
\times\,\frac{\partial n_s}{\partial y^i}\,
\frac{\partial n_r}{\partial y^j}
+\sum^n_{s=1}\sum^n_{r=1}\left(\,\shave{\sum^n_{q=1}}
\nu\,p_q\,\frac{\partial}{\partial p_q}\!
\left(\frac{1}{\Omega}\,\frac{\partial H}{\partial x^r}
-Q_r\!\right)\!\left(\frac{1}{\Omega}\,\frac{\partial H}
{\partial x^s}-Q_s\!\right)\,-\right.\\
\vspace{1ex}
\left.\vphantom{\shave{\sum^n_{q=1}}}
-\,\nu\,\frac{\partial}{\partial x^s}\!
\left(\frac{1}{\Omega}\,\frac{\partial H}{\partial x^r}
-Q_r\!\right)\right)\tau^s_i\,\tau^r_j.
\endgather
$$
Term with product of partial derivatives $\partial n_s/\partial y^i$
and $\partial n_r/\partial y^j$ in the above expression can be replaced
by dots. Indeed, one can easily check up that coefficient of such product
of derivatives is symmetric in indices $s$ and $r$:
$$
\gather
\sum^n_{q=1}
\frac{\nu^3\,p_q}{\Omega}\,\frac{\partial}{\partial p_q}\!
\left(\frac{1}{\Omega}\,\frac{\partial H}{\partial p_r}\right)
\frac{\partial H}{\partial p_s}+\frac{\nu^3}{\Omega^2}\,  
\frac{\partial\Omega}{\partial p_s}\,\frac{\partial H}
{\partial p_r}=
-\sum^n_{q=1}\frac{\nu^3\,p_q}{\Omega^3}\,\frac{\partial\Omega}
{\partial p_q}\,\frac{\partial H}{\partial p_r}\,\frac{\partial H}
{\partial p_s}\,+\\
\vspace{1ex}
+\sum^n_{q=1}
\frac{\nu^3\,p_q}{\Omega^2}\,\frac{\partial^2H}{\partial p_q\,
\partial p_r}\,\frac{\partial H}{\partial p_s}
+\sum^n_{q=1}
\frac{\nu^3\,p_q}{\Omega^2}\,\frac{\partial^2H}{\partial p_q\,
\partial p_s}\,\frac{\partial H}{\partial p_r}
+\frac{\nu^3}{\Omega^2}\,\frac{\partial H}{\partial p_s}\,
\frac{\partial H}{\partial p_r}.
\endgather
$$
Now let's study the term with $\tau^s_i\,\tau^r_j$. For the coefficient
in this term we derive
$$
\gather
\sum^n_{q=1}\nu\,p_q\,\frac{\partial}{\partial p_q}\!
\left(\frac{1}{\Omega}\,\frac{\partial H}{\partial x^r}
-Q_r\!\right)\!\left(\frac{1}{\Omega}\,\frac{\partial H}
{\partial x^s}-Q_s\!\right)
-\nu\,\frac{\partial}{\partial x^s}\!
\left(\frac{1}{\Omega}\,\frac{\partial H}{\partial x^r}
-Q_r\!\right)=\\
\displaybreak
=\sum^n_{q=1}\nu\,p_q\left(-\frac{1}{\Omega^2}\,
\frac{\partial\Omega}{\partial p_q}\,\frac{\partial H}
{\partial x^r}+\frac{1}{\Omega}\,\frac{\partial^2 H}
{\partial p_q\,\partial x^r}-\frac{\partial Q_r}
{\partial p_q}\right)\!\left(\frac{1}{\Omega}\,
\frac{\partial H}{\partial x^s}-Q_s\!\right)-\\
\vspace{1ex}
-\,\nu\left(-\frac{1}{\Omega^2}\,\frac{\partial\Omega}
{\partial x^s}\,\frac{\partial H}{\partial x^r}+
\frac{1}{\Omega}\,\frac{\partial^2 H}{\partial x^s\,
\partial x^r}-\frac{\partial Q_r}{\partial x^s}\right)
=\ldots+\sum^n_{q=1}\frac{\nu\,p_q}{\Omega^2}\,
\frac{\partial^2 H}{\partial p_q\,\partial x^r}\,
\frac{\partial H}{\partial x^s}\,-\\
\vspace{1ex}
-\sum^n_{q=1}\frac{\nu\,p_q}{\Omega}\,\frac{\partial Q_r}
{\partial p_q}\,\frac{\partial H}{\partial x^s}
+\sum^n_{q=1}\frac{\nu\,p_q}{\Omega^2}\,\frac{\partial H}
{\partial p_q}\,\frac{\partial H}{\partial x^r}\,Q_s
+\sum^n_{q=1}\sum^n_{k=1}\frac{\nu\,p_q\,p_k}{\Omega^2}
\frac{\partial^2H}{\partial p_q\,\partial p_k}\,
\frac{\partial H}{\partial x^r}\,Q_s\,-\\
\vspace{1ex}
-\sum^n_{q=1}\frac{\nu\,p_q}{\Omega}\,\frac{\partial^2H}
{\partial p_q\,\partial x^r}\,Q_s+\sum^n_{q=1}\nu\,p_q\,
\frac{\partial Q_r}{\partial p_q}\,Q_s+\sum^n_{q=1}
\frac{\nu\,p_q}{\Omega^2}\,\frac{\partial^2H}
{\partial p_q\,\partial x^s}\,\frac{\partial H}
{\partial x^r}-\frac{\nu}{\Omega}\,\frac{\partial^2H}
{\partial x^r\,\partial x^s}\,+\\
\vspace{1ex}
+\,\nu\,\frac{\partial Q_r}{\partial x^s}=\dots
-\sum^n_{q=1}\frac{\nu\,p_q}{\Omega}\,\frac{\partial Q_r}
{\partial p_q}\,\frac{\partial H}{\partial x^s}
+\sum^n_{q=1}\frac{\nu\,p_q}{\Omega^2}\,\frac{\partial H}
{\partial p_q}\,\frac{\partial H}{\partial x^r}\,Q_s
+\sum^n_{q=1}\sum^n_{k=1}\frac{\nu\,p_q\,p_k}{\Omega^2}
\,\times\\
\vspace{1ex}
\times\,\frac{\partial^2H}{\partial p_q\,\partial p_k}\,
\frac{\partial H}{\partial x^r}\,Q_s
-\sum^n_{q=1}\frac{\nu\,p_q}{\Omega}\,\frac{\partial^2H}
{\partial p_q\,\partial x^r}\,Q_s+\sum^n_{q=1}\nu\,p_q\,
\frac{\partial Q_r}{\partial p_q}\,Q_s
+\nu\,\frac{\partial Q_r}{\partial x^s}.
\endgather
$$
In the above calculations we used formula \thetag{7.7} for $\Omega$.
As a result for partial derivative \thetag{7.17} we have derived the
following expression:
$$
\gather
\frac{\partial^2\nu}{\partial y^i\,\partial y^j}=\dots
+\sum^n_{s=1}\sum^n_{r=1}\left(\frac{2\,\nu^2}{\Omega}
\left(\frac{1}{\Omega}\,\frac{\partial H}{\partial x^s}
-Q_s\!\right)\frac{\partial H}
{\partial p_r}-\nu^2\,\frac{\partial}{\partial x^s}\!
\left(\frac{1}{\Omega}\,\frac{\partial H}{\partial p_r}\right)
+\shave{\sum^n_{q=1}}\nu^2\,\times
\vphantom{\shave{\sum^n_{q=1}}}\right.\\
\vspace{1ex}
\left.\times\,p_q\,\frac{\partial}
{\partial p_q}\!\left(\frac{1}{\Omega}\,\frac{\partial H}
{\partial p_r}\right)\!\left(\frac{1}{\Omega}\,\frac{\partial H}
{\partial x^s}-Q_s\!\right)\right)
\frac{\partial n_r}{\partial y^j}\,\tau^s_i
+\sum^n_{s=1}\sum^n_{r=1}\left(\frac{\nu^2}{\Omega}\,
\frac{\partial H}{\partial p_s}\,\left(\frac{1}{\Omega}\,
\frac{\partial H}{\partial x^r}-Q_r\!\right)-
\vphantom{\shave{\sum^n_{q=1}}}\right.\\
\vspace{1ex}
\left.-\,\nu^2\,\frac{\partial}{\partial p_s}\!
\left(\frac{1}{\Omega}\,\frac{\partial H}{\partial x^r}
-Q_r\!\right)
+\shave{\sum^n_{q=1}}\frac{\nu^2\,p_q}{\Omega}\,
\frac{\partial}{\partial p_q}\!\left(\frac{1}{\Omega}\,
\frac{\partial H}{\partial x^r}-Q_r\!\right)
\frac{\partial H}{\partial p_s}\right)
\frac{\partial n_s}{\partial y^i}\,\tau^r_j\,+\\
\vspace{1ex}
+\sum^n_{s=1}\sum^n_{r=1}\left(-\shave{\sum^n_{q=1}}
\frac{\nu\,p_q}{\Omega}\,\frac{\partial Q_r}{\partial p_q}
\,\frac{\partial H}{\partial x^s}+\shave{\sum^n_{q=1}}
\frac{\nu\,p_q}{\Omega^2}\,\frac{\partial H}{\partial p_q}
\,\frac{\partial H}{\partial x^r}\,Q_s
+\sum^n_{q=1}\sum^n_{k=1}\frac{\nu\,p_q\,p_k}{\Omega^2}\,
\times\right.\\
\vspace{1ex}
\left.\times\,\frac{\partial^2H}{\partial p_q\,\partial p_k}\,
\frac{\partial H}{\partial x^r}\,Q_s
-\sum^n_{q=1}\frac{\nu\,p_q}{\Omega}\,\frac{\partial^2H}
{\partial p_q\,\partial x^r}\,Q_s
+\shave{\sum^n_{q=1}}\nu\,p_q\,
\frac{\partial Q_r}{\partial p_q}\,Q_s
+\nu\,\frac{\partial Q_r}{\partial x^s}
\right)\tau^s_i\,\tau^r_j.
\endgather
$$
Next step in transforming the above expression is based on the
following equality:
$$
\sum^n_{s=1}\sum^n_{r=1}\alpha^r_s\,\frac{\partial n_r}
{\partial y^j}\,\tau^s_i+\sum^n_{s=1}\sum^n_{r=1}\beta^s_r\,
\frac{\partial n_s}{\partial y^i}\,\tau^r_j=\sum^n_{s=1}
\sum^n_{r=1}(\alpha^r_s-\beta^r_s)\,\frac{\partial n_r}
{\partial y^j}\,\tau^s_i+\ldots\ .
$$
Here, as we already used above, we denoted by dots terms symmetric
in indices $i$ and $j$. Further for our particular $\alpha_{rs}$ and
$\beta_{rs}$ we derive
$$
\gather
\alpha^r_s-\beta^r_s=\frac{\nu^2}{\Omega}
\left(\frac{1}{\Omega}\,\frac{\partial H}{\partial x^s}
-Q_s\!\right)\frac{\partial H}{\partial p_r}-\nu^2\,
\frac{\partial}{\partial x^s}\!\left(\frac{1}{\Omega}\,
\frac{\partial H}{\partial p_r}\right)
+\sum^n_{q=1}\nu^2\,p_q\,\times\\
\displaybreak
\times\,\frac{\partial}
{\partial p_q}\!\left(\frac{1}{\Omega}\,\frac{\partial H}
{\partial p_r}\right)\!\left(\frac{1}{\Omega}\,\frac{\partial H}
{\partial x^s}-Q_s\!\right)+\nu^2\,\frac{\partial}{\partial p_r}\!
\left(\frac{1}{\Omega}\,\frac{\partial H}{\partial x^s}
-Q_s\!\right)-\sum^n_{q=1}\frac{\nu^2\,p_q}
{\Omega}\,\times\\
\vspace{1ex}
\times\,\frac{\partial}{\partial p_q}\!\left(\frac{1}{\Omega}\,
\frac{\partial H}{\partial x^s}-Q_s\!\right)\frac{\partial H}
{\partial p_r}=\sum^n_{q=1}\sum^n_{k=1}\frac{\nu^2}{\Omega^2}
\,p_q\,p_k\,\frac{\partial^2H}{\partial p_q\,\partial p_k}\,
\frac{\partial H}{\partial p_r}\,Q_s\,-\\
\vspace{1ex}
-\sum^n_{q=1}\frac{\nu^2}{\Omega}\,p_q\,\frac{\partial^2H}
{\partial p_q\,\partial p_r}\,Q_s-\nu^2\,\frac{\partial Q_s}
{\partial p_r}+\sum^n_{q=1}\frac{\nu^2}{\Omega}\,p_q\,
\frac{\partial H}{\partial p_r}\,\frac{\partial Q_s}{\partial p_q}.
\endgather
$$
Summarizing all above calculations, for partial derivatives
\thetag{7.17} we obtain
$$
\hskip -0.5em
\gathered
\frac{\partial^2\nu}{\partial y^i\,\partial y^j}=\dots
+\sum^n_{s=1}\sum^n_{r=1}\left(\,\shave{\sum^n_{q=1}}
\shave{\sum^n_{k=1}}\frac{\nu^2}{\Omega^2}\,p_q\,p_k\,
\frac{\partial^2H}{\partial p_q\,\partial p_k}\,
\frac{\partial H}{\partial p_r}\,Q_s\,-\right.\\
\vspace{1ex}
\left.-\shave{\sum^n_{q=1}}\frac{\nu^2}{\Omega}\,p_q\,\frac{\partial^2H}
{\partial p_q\,\partial p_r}\,Q_s-\nu^2\,\frac{\partial Q_s}
{\partial p_r}+\shave{\sum^n_{q=1}}\frac{\nu^2}{\Omega}\,p_q\,
\frac{\partial H}{\partial p_r}\,\frac{\partial Q_s}{\partial p_q}
\right)\frac{\partial n_r}{\partial y^j}\,\tau^s_i\,+\\
\vspace{1ex}
+\sum^n_{s=1}\sum^n_{r=1}\left(-\shave{\sum^n_{q=1}}
\frac{\nu\,p_q}{\Omega}\,\frac{\partial Q_r}{\partial p_q}
\,\frac{\partial H}{\partial x^s}+\shave{\sum^n_{q=1}}
\frac{\nu\,p_q}{\Omega^2}\,\frac{\partial H}{\partial p_q}
\,\frac{\partial H}{\partial x^r}\,Q_s
+\shave{\sum^n_{q=1}}\shave{\sum^n_{k=1}}\frac{\nu\,p_q\,p_k}{\Omega^2}\,
\times\right.\\
\vspace{1ex}
\left.\times\,\frac{\partial^2H}{\partial p_q\,\partial p_k}\,
\frac{\partial H}{\partial x^r}\,Q_s
-\shave{\sum^n_{q=1}}\frac{\nu\,p_q}{\Omega}\,\frac{\partial^2H}
{\partial p_q\,\partial x^r}\,Q_s
+\shave{\sum^n_{q=1}}\nu\,p_q\,
\frac{\partial Q_r}{\partial p_q}\,Q_s
+\nu\,\frac{\partial Q_r}{\partial x^s}
\right)\tau^s_i\,\tau^r_j.
\endgathered\hskip -1em
\tag9.1
$$
In a similar way for partial derivative \thetag{7.18} one can get
analogous expression:
$$
\hskip -0.5em
\gathered
\frac{\partial^2\nu}{\partial y^j\,\partial y^i}=\dots
+\sum^n_{s=1}\sum^n_{r=1}\left(\shave{\sum^n_{q=1}}
\shave{\sum^n_{k=1}}\frac{\nu^2}{\Omega^2}\,p_q\,p_k\,
\frac{\partial^2H}{\partial p_q\,\partial p_k}\,
\frac{\partial H}{\partial p_s}\,Q_r\,-\right.\\
\vspace{1ex}
\left.-\shave{\sum^n_{q=1}}\frac{\nu^2}{\Omega}\,p_q\,\frac{\partial^2H}
{\partial p_q\,\partial p_s}\,Q_r-\nu^2\,\frac{\partial Q_r}
{\partial p_s}+\shave{\sum^n_{q=1}}\frac{\nu^2}{\Omega}\,p_q\,
\frac{\partial H}{\partial p_s}\,\frac{\partial Q_r}{\partial p_q}
\right)\frac{\partial n_s}{\partial y^i}\,\tau^r_j\,+\\
\vspace{1ex}
+\sum^n_{s=1}\sum^n_{r=1}\left(-\shave{\sum^n_{q=1}}
\frac{\nu\,p_q}{\Omega}\,\frac{\partial Q_s}{\partial p_q}
\,\frac{\partial H}{\partial x^r}+\shave{\sum^n_{q=1}}
\frac{\nu\,p_q}{\Omega^2}\,\frac{\partial H}{\partial p_q}
\,\frac{\partial H}{\partial x^s}\,Q_r
+\shave{\sum^n_{q=1}}\shave{\sum^n_{k=1}}\frac{\nu\,p_q\,p_k}
{\Omega^2}\,\times\right.\\
\vspace{1ex}
\left.\times\,\frac{\partial^2H}{\partial p_q\,\partial p_k}\,
\frac{\partial H}{\partial x^s}\,Q_r
-\shave{\sum^n_{q=1}}\frac{\nu\,p_q}{\Omega}\,\frac{\partial^2H}
{\partial p_q\,\partial x^s}\,Q_r
+\shave{\sum^n_{q=1}}\nu\,p_q\,
\frac{\partial Q_s}{\partial p_q}\,Q_r
+\nu\,\frac{\partial Q_s}{\partial x^r}
\right)\tau^s_i\,\tau^r_j.
\endgathered\hskip -1em
\tag9.2
$$
Now we are able to write compatibility equations \thetag{7.19} in
explicit form. It is sufficient to equate partial derivatives
\thetag{9.1} and \thetag{9.2} to each other. However, we shall not
do this right now because we need some additional information in
order to treat arising compatibility equation properly. We should
replace partial derivatives in \thetag{9.1} and \thetag{9.2} by
covariant derivatives, and we should understand some facts from
theory of hypersurfaces in a manifold equipped with the only
geometric structure given by Lagrange function $L$.
\head
10. Differentiation of extended tensor fields.
\endhead
    Let's consider smooth extended tensor fields in their
$\bold p$-representation as given by definition~5.1. They form
a graded ring with respect to standard operations of addition
and tensor product. We denote it as follows:
$$
\hskip -2em
\bold T(M)=\bigoplus^\infty_{r=0}\bigoplus^\infty_{s=0}T^r_s(M).
\tag10.1
$$
Graded ring \thetag{10.1} is equipped with operation of contraction,
which is also standard. Moreover, \thetag{10.1} possesses the structure
of algebra over the ring of smooth functions in $T^*\!M$. For this
reason it is called {\it extended algebra of tensor fields}.\par
    In extended algebra of tensor fields \thetag{10.1} one can define
canonical covariant differentiation, which is called {\it vertical
gradient} or {\it momentum gradient}:
$$
\hskip -2em
\tilde\nabla\!:\ T^r_s(M)\to T^{r+1}_s(M).
\tag10.2
$$
In local chart momentum gradient \thetag{10.2} is determined by the
following formula:
$$
\hskip -2em
\tilde\nabla^qX^{i_1\ldots\,i_r}_{j_1\ldots\,j_s}=
\frac{\partial X^{i_1\ldots\,i_r}_{j_1\ldots\,j_s}}
{\partial p_q}.
\tag10.3
$$
To define another covariant differentiation in $\bold T(M)$ one need
some additional geometric structure in $M$. This is so called {\it
extended affine connection}. We shall not discuss the nature
of this structure (see thesis \cite{6}). Note only that in local
chart it is given by its components $\Gamma^k_{ij}$, which obey standard
transformation rule:
$$
\hskip -2em
\Gamma^k_{ij}=\sum^n_{m=1}\sum^n_{a=1}\sum^n_{c=1} S^k_m\,T^a_i
\,T^c_j\,\tilde\Gamma^m_{ac}+\sum^n_{m=1} S^k_m\,\frac{\partial
T^m_i}{\partial x^j}.
\tag10.4
$$
Here $S^i_j$ and $T^i_j$ are components of transition matrices
$S$ and $T$ binding coordinates $x^1,\,\ldots,\,x^n$ and
$\tx^1,\,\ldots,\,\tx^n$ in two overlapping local charts of $M$:
$$
\xalignat 2
&\hskip -2em
S^i_j=\frac{\partial x^i}{\partial\tx^j},
&T^i_j=\frac{\partial\tx^i}{\partial x^j}.
\tag10.5
\endxalignat
$$
Unlike components of standard affine connection, components of
extended affine connection $\Gamma$ in $\bold p$-representation
depend not only on $x^1,\,\ldots,\,x^n$, but also on components
$p_1,\,\ldots,\,p_n$ of momentum covector $\bold p$:
$$
\hskip -2em
\Gamma^k_{ij}=\Gamma^k_{ij}(x^1,\ldots,x^n,p_1,\ldots,p_n).
\tag10.6
$$
Thus, if $M$ possesses extended affine connection $\Gamma$, one
can define {\it horizontal gradient} or {\it spatial gradient}
$\nabla$ in $\bold T(M)$. In local chart it is expressed by 
formula
$$
\hskip -2em
\aligned
&\nabla_{\!q}X^{i_1\ldots\,i_r}_{j_1\ldots\,j_s}=\frac{\partial
X^{i_1\ldots\,i_r}_{j_1\ldots\,j_s}}{\partial x^q}
+\sum^n_{a=1}\sum^n_{b=1}p_a\,\Gamma^a_{qb}\,\frac{\partial
X^{i_1\ldots\,i_r}_{j_1\ldots\,j_s}}{\partial p_b}\,+\\
&+\sum^r_{k=1}\sum^n_{a_k=1}\!\Gamma^{i_k}_{q\,a_k}\,X^{i_1\ldots\,
a_k\ldots\,i_r}_{j_1\ldots\,\ldots\,\ldots\,j_s}
-\sum^s_{k=1}\sum^n_{b_k=1}\!\Gamma^{b_k}_{q\,j_k}\,
X^{i_1\ldots\,\ldots\,\ldots\,i_r}_{j_1\ldots\,b_k\ldots\,j_s}.
\endaligned
\tag10.7
$$
Horizontal gradient $\nabla$ increments by $1$ the number of lower
indices of tensor field:
$$
\hskip -2em
\nabla\!:\ T^r_s(M)\to T^r_{s+1}(M).
\tag10.8
$$\par
    Using Legendre map, one can transform extended tensor fields to
$\bold v$-representa\-tion. In $\bold v$-representation vertical
gradient $\tilde\nabla$ is called {\it velocity gradient}. It is
a~map
$$
\hskip -2em
\tilde\nabla\!:\ T^r_s(M)\to T^r_{s+1}(M)
\tag10.9
$$
similar to \thetag{10.8}. In local chart this map \thetag{10.9}
is given by formula
$$
\hskip -2em
\tilde\nabla_{\!q}X^{i_1\ldots\,i_r}_{j_1\ldots\,j_s}=\frac{\partial
X^{i_1\ldots\,i_r}_{j_1\ldots\,j_s}}{\partial v^q}.
\tag10.10
$$
Two versions of vertical gradient \thetag{10.3} and \thetag{10.10}
can be related to each other by means of formula using components
of matrix \thetag{6.8}:
$$
\hskip -2em
\tilde\nabla_{\!q}\!\left(X^{i_1\ldots\,i_r}_{j_1\ldots\,j_s}\compos
\lambda\right)=\sum^n_{k=1}g_{qk}\,\tilde\nabla^k
X^{i_1\ldots\,i_r}_{j_1\ldots\,j_s}\compos\lambda.
\tag10.11
$$
Remember that inverse Legendre map in local chart is given by
explicit formula
$$
\hskip -2em
v^i=\frac{\partial H}{\partial p_i}.
\tag10.12
$$
Therefore we can introduce matrix $\bold g^{-1}$ with the
following components:
$$
\hskip -2em
g^{ij}=\frac{\partial^2\!H}{\partial p_i\,\partial p_j}.
\tag10.13
$$
Matrices \thetag{6.8} and \thetag{10.13} are inverse to each other
when taken in the same representation, i\.\,e\. when both brought
to $\bold p$-representation or when both brought to
$\bold v$-representation. Matrix \thetag{10.13} is used in formula,
which is similar to \thetag{10.11}:
$$
\hskip -2em
\tilde\nabla^q\!\left(X^{i_1\ldots\,i_r}_{j_1\ldots\,j_s}\compos
\lambda^{-1}\right)=\sum^n_{k=1}g^{qk}\,\tilde\nabla_{\!k}
X^{i_1\ldots\,i_r}_{j_1\ldots\,j_s}\compos\lambda^{-1}.
\tag10.14
$$
Looking at formulas \thetag{10.11} and \thetag{10.14}, we see that
matrices \thetag{6.8} and \thetag{10.13} here do part of work that
metric tensor does in Riemannian geometry.\par
    In order to define horizontal gradient for extended tensor fields
in $\bold v$-representa\-tion we should transform connection components
\thetag{10.6} by means of Legendre map:
$$
\hskip -2em
\Gamma^k_{ij}\to\Gamma^k_{ij}\compos\lambda=\Gamma^k_{ij}(x^1,
\ldots,x^n,v^1,\ldots,v^n).
\tag10.15
$$
Then, using \thetag{10.15}, we can define horizontal gradient by
a formula in local chart:
$$
\hskip -2em
\aligned
&\nabla_{\!q}X^{i_1\ldots\,i_r}_{j_1\ldots\,j_s}=\frac{\partial
X^{i_1\ldots\,i_r}_{j_1\ldots\,j_s}}{\partial x^q}
-\sum^n_{a=1}\sum^n_{b=1}v^a\,\Gamma^b_{qa}\,\frac{\partial
X^{i_1\ldots\,i_r}_{j_1\ldots\,j_s}}{\partial v^b}\,+\\
&+\sum^r_{k=1}\sum^n_{a_k=1}\!\Gamma^{i_k}_{q\,a_k}\,X^{i_1\ldots\,
a_k\ldots\,i_r}_{j_1\ldots\,\ldots\,\ldots\,j_s}
-\sum^s_{k=1}\sum^n_{b_k=1}\!\Gamma^{b_k}_{q\,j_k}
X^{i_1\ldots\,\ldots\,\ldots\,i_r}_{j_1\ldots\,b_k\ldots\,j_s}.
\endaligned
\tag10.16
$$
In spite of transformation rule \thetag{10.15} for connection
components, in general, \thetag{10.7} and \thetag{10.16} are
not different representations of the same tensor field. We have
the following equality binding these two horizontal gradients:
$$
\hskip -2em
\nabla_{\!q}\!\left(X^{i_1\ldots\,i_r}_{j_1\ldots\,j_s}\compos
\lambda\right)-\nabla_{\!q}X^{i_1\ldots\,i_r}_{j_1\ldots\,j_s}
\compos\lambda=\sum^n_{s=1}\nabla_{\!q}\!\tilde\nabla_{\!s}L\,
\cdot\tilde\nabla^s\!X^{i_1\ldots\,i_r}_{j_1\ldots\,j_s}\compos
\lambda.
\tag10.17
$$
However, there is special case, when \thetag{10.7} and \thetag{10.16}
do coincide. This is when
$$
\hskip -2em
\nabla_{\!q}\!\tilde\nabla_{\!s}L=0.
\tag10.18
$$
In this case \thetag{10.17} turns to equality, which is similar to
the equality \thetag{10.11}:
$$
\hskip -2em
\nabla_{\!q}\!\left(X^{i_1\ldots\,i_r}_{j_1\ldots\,j_s}\compos
\lambda\right)=\nabla_{\!q}X^{i_1\ldots\,i_r}_{j_1\ldots\,j_s}
\compos\lambda.
\tag10.19
$$
\definition{Definition 10.1} Extended connection $\Gamma$ is
called {\it concordant} with Lagrange function $L$ if the
equality \thetag{10.18} holds.
\enddefinition
For concordant connections due to \thetag{10.19} horizontal 
gradient can be calculated either in $\bold p$ or in
$\bold v$-representation yielding the same result, but in
different variables. Note that the equality \thetag{10.18} can
be replaced by equivalent equality for Hamilton function $H$
in $\bold p$-representation. It looks like
$$
\hskip -2em
\nabla_{\!q}\tilde\nabla^sH=0.
\tag10.20
$$
This equality \thetag{10.20} can be derived from \thetag{10.18}
by direct calculations.
\head
11. Differentiation along lines and hypersurfaces.
\endhead
    Let $p=p(t)$ be a smooth parametric curve in our manifold $M$. In
local chart with coordinates $x^1,\,\ldots,\,x^n$ it is represented by
functions
$$
\hskip -2em
\cases
x^1=x^1(t),\\
.\ .\ .\ .\ .\ .\ .\ .\\
x^n=x^n(t).
\endcases
\tag11.1
$$
Suppose that at each point $p$ of this curve \thetag{6.1} some tensor
of the type $(r,s)$ is given, i\.\,e\. we have tensor-valued function
$\bold X=\bold X(t)$. In local chart this tensor function is expressed
by its components $X^{i_1\ldots\,i_r}_{j_1\ldots\,j_s}=
X^{i_1\ldots\,i_r}_{j_1\ldots\,j_s}(t)$. If we had standard affine
connection in $M$, we could define covariant derivative of tensor
function $\bold X(t)$ with respect to parameter $t$ along the curve.
This is another tensor function $\bold Y=\nabla_{\!t}\bold X$ with
components $Y^{i_1\ldots\,i_r}_{j_1\ldots\,j_s}=
\nabla_{\!t}X^{i_1\ldots\,i_r}_{j_1\ldots\,j_s}$ given by the following
formula:
$$
\hskip -2em
\gathered
\nabla_{\!t}X^{i_1\ldots\,i_r}_{j_1\ldots\,j_s}=
\frac{dX^{i_1\ldots\,i_r}_{j_1\ldots\,j_s}}{dt}
+\sum^r_{k=1}\sum^n_{m=1}\sum^n_{a_k=1}\dot x^m\,\Gamma^{i_k}_{m\,a_k}
\,X^{i_1\ldots\,a_k\ldots\,i_r}_{j_1\ldots\,\ldots\,\ldots\,j_s}\,-\\
\vspace{1ex}
-\sum^s_{k=1}\sum^n_{m=1}\sum^n_{b_k=1}\dot x^m\,\Gamma^{b_k}_{m\,j_k}
X^{i_1\ldots\,\ldots\,\ldots\,i_r}_{j_1\ldots\,b_k\ldots\,j_s}.
\endgathered
\tag11.2
$$
However, with extended connection $\Gamma$ we need come additional
data in order to apply formula \thetag{11.2} to tensor function
$\bold X=\bold X(t)$. Indeed, looking at \thetag{10.6} and
\thetag{10.15}, we see that in $\bold p$-representation we need
functions
$$
\hskip -2em
p_1(t),\ \ldots,\ p_n(t),
\tag11.3
$$
while in $\bold v$-representation we need another set of functions:
$$
\hskip -2em
v^1(t),\ \ldots,\ v^n(t).
\tag11.4
$$
Taken by themselves, functions \thetag{11.3} and \thetag{11.4} are
components of covector-function $\bold p=\bold p(t)$ and
vector-function $\bold v=\bold v(t)$. But combined with functions
\thetag{11.1}, they define lift of initial curve $p=p(t)$ from $M$
to cotangent bundle $T^*\!M$ and to tangent bundle $TM$ respectively.
Once such lift is given, we can use formula \thetag{11.2}. Covariant
derivative given by this formula then is called {\it covariant
derivative with respect to parameter $t$ along the curve $p=p(t)$
due to its lift $q=q(t)$}.\par
    Now let's discuss the question of how tensor function
$\bold X=\bold X(t)$ on curve could be defined. Surely, it
can be given immediately as it is. However, often tensor function
$\bold X=\bold X(t)$ comes to curve from outer space, i\.\,e\.
from manifold $M$. For example, if in $M$ or at least in some
neighborhood of our curve some standard (not extended) tensor
field $\bold X$ is given, we can restrict it to the curve $p=p(t)$
thus obtaining tensor function $\bold X=\bold X(t)$. In the case of
extended tensor field $\bold X$ we cannot restrict it to the curve
$p=p(t)$ immediately. We should first choose a lift of this curve
$q=q(t)$ in $TM$ or in $T^*\!M$, then we could restrict $\bold X$
to the lifted curve. Now suppose that tensor function $\bold X=
\bold X(t)$ on the curve $p=p(t)$ is obtained from extended tensor
field $\bold X$ in this way. Then covariant derivative \thetag{11.2}
for $\bold X(t)$ is determined. By direct calculations we can prove
that this covariant derivative can be expressed through covariant
derivatives \thetag{10.7} and \thetag{10.3} in $\bold p$-representation:
$$
\hskip -2em
\nabla_{\!t}X^{i_1\ldots\,i_r}_{j_1\ldots\,j_s}=\sum^n_{k=1}\nabla_{\!k}
X^{i_1\ldots\,i_r}_{j_1\ldots\,j_s}\cdot\dot x^k+\sum^n_{k=1}
\tilde\nabla^k\!X^{i_1\ldots\,i_r}_{j_1\ldots\,j_s}\cdot \nabla_{\!t}p_k.
\tag11.5
$$
In $\bold v$-representation we have similar expression
$$
\hskip -2em
\nabla_{\!t}X^{i_1\ldots\,i_r}_{j_1\ldots\,j_s}=\sum^n_{k=1}\nabla_{\!k}
X^{i_1\ldots\,i_r}_{j_1\ldots\,j_s}\cdot\dot x^k+\sum^n_{k=1}
\tilde\nabla_{\!k}X^{i_1\ldots\,i_r}_{j_1\ldots\,j_s}\cdot\nabla_{\!t}v^k
\tag11.6
$$
that uses covariant derivatives \thetag{10.16} and \thetag{10.10}.
Here $\nabla_{\!t}v^k$ and $\nabla_{\!t}p_k$ are components of
covariant derivatives for vectorial and covectorial functions
with components \thetag{11.4} and \thetag{11.3} that determine
lift of curve $p=p(t)$ to $TM$ and $T^*\!M$ respectively. From
general point of view, formulas \thetag{11.5} and \thetag{11.6}
are nothing, but the rule of differentiating composite functions
written in terms of covariant derivatives. They are obtained by
direct calculations.\par
\parshape 26 0pt 360pt 0pt 360pt
0pt 360pt 0pt 360pt 0pt 360pt 0pt 360pt 0pt 360pt
0pt 360pt 140pt 220pt
140pt 220pt 140pt 220pt 140pt 220pt 140pt 220pt 140pt 220pt
140pt 220pt 140pt 220pt 140pt 220pt 140pt 220pt 140pt 220pt
140pt 220pt 140pt 220pt 140pt 220pt 140pt 220pt 140pt 220pt
140pt 220pt 0pt 360pt
    Now let $S$ be a smooth hypersurface in $M$. Denote by $y^1,
\,\ldots,\,y^{n-1}$ coordinates of point in some local chart of
$S$. Then the following smooth functions
$$
\hskip -2em
\cases
x^1=x^1(y^1,\ldots,y^{n-1}),\\
.\ .\ .\ .\ .\ .\ .\ .\ .\ .\ .\ .\ .\ .\ .\ .\ .\\
x^n=x^n(y^1,\ldots,y^{n-1})
\endcases
\tag11.7
$$
represent hypersurface $S$ in local chart of $M$. Partial
derivatives of functions \thetag{11.7} determine tangent
vectors $\boldsymbol\tau_1,\,\ldots,\,\boldsymbol\tau_{n-1}$
to $S$ as it is done in \thetag{2.5}. Taking one variable
$y^i$ in \thetag{11.7} and fixing all others, we can consider
\thetag{11.7} as $n-1$ families of functions of one variable,
\vadjust{\vskip 1pt\hbox to 0pt{\kern 10pt\hbox{\special{em:graph
pst-14b.gif}}\hss}\vskip -1pt}each corresponding to some family
of curves on $S$. These are
coordinate curves forming coordinate network on $S$ (see Fig.~11.1).
One can lift coordinate curves to cotangent bundle $T^*\!M$ by means
of normal covector $\bold n=\bold n(p)$. However, we shall use another
lift given by momentum covector (see initial data \thetag{7.4}):
$$
\hskip -2em
\bold p=\nu(p)\cdot\bold n(p).
\tag11.8
$$
Then applying Legendre transformation \thetag{1.5} to covector
\thetag{11.8}, we obtain vector
$$
\hskip -2em
\bold v=\bold v(p)=\tilde\nabla L
\tag11.9
$$
at each point $p\in S$. Vector $\bold v(p)$ is transversal to $S$
as shown on Fig.~11.1. This follows from regularity of Lagrange
function (see definition~1.1). Now, using lifts defined by \thetag{11.8}
and \thetag{11.9}, we can apply formula \thetag{11.2} to coordinate
lines of hypersurface $S$.\par
    Remember that $t=y^i$ is a parameter of $i$-th coordinate line,
while $\boldsymbol\tau_i$ is tangent vector to this line corresponding
to parameter $y^i$. Covariant derivative $\nabla_{\!t}$ with respect
to parameter $t=y^i$ along $i$-th coordinate line by tradition is
denoted by $\nabla_{\!\boldsymbol\tau_i}$. This derivative can be
applied to any smooth tensorial function defined at the points of
hypersurface $S$. Let's apply it to vector-function $\boldsymbol
\tau_j$. As a result ge get another vector-function $\nabla_{\!
\boldsymbol\tau_i}\boldsymbol\tau_j$ on $S$. Similarly one can
consider vector-function $\nabla_{\!\boldsymbol\tau_i}\bold v$
and covector-function $\nabla_{\!\boldsymbol\tau_i}\bold p$ on
$S$. The latter one appears to be most important for us. It will
be studied in section~13 below.
\head
12. Covariant form of compatibility equations.
\endhead
   In section~7 we have derived Pfaff equations \thetag{7.13} for
scalar function $\nu$ on hypersurface $S$. Then in section~9 we
studied compatibility condition for these Pfaff equations. Explicit
form of compatibility condition could be obtained by equating
partial derivatives \thetag{9.1} and \thetag{9.2}. But as a result
we would obtain huge equality difficult to observe. This means that
formulas \thetag{9.1} and \thetag{9.2} require some preliminary
transformations in order to simplify further analysis of compatibility
condition they lead to.\par
    First of all we are going to replace partial derivatives $\partial
n_r/\partial y^j$ and $\partial n_s/\partial y^i$ by covariant derivatives
$\nabla_{\!\boldsymbol\tau_j}p_r$ and $\nabla_{\!\boldsymbol\tau_i}p_s$,
assuming that we have some symmetric extended affine connection in $M$.
Let's turn back to the equalities \thetag{7.12} and \thetag{7.13},
let's substitute $\partial\nu/\partial y^i$ from \thetag{7.13} into
\thetag{7.12} and remember formula \thetag{10.12}:
$$
\pagebreak
\hskip -2em
\frac{\partial p_s}{\partial y^i}
=\sum^n_{r=1}\nu\left(\delta^r_s-\frac{v^r\,p_s}{\Omega}\right)
\frac{\partial n_r}{\partial y^i}
-\sum^n_{r=1}\left(\frac{1}{\Omega}\frac{\partial H}
{\partial x^r}-Q_r\!\right)\tau^r_i\,p_s.
\tag12.1
$$
Denote by $P^r_s$ coefficient in front of partial derivative
$\partial n_r/\partial y^i$ in formula \thetag{12.1}:
$$
\hskip -2em
P^r_s=\delta^r_s-\frac{v^r\,p_s}{\Omega}.
\tag12.2
$$
Formula \thetag{12.2} for $P^r_s$ can be rewritten in two equivalent
forms:
$$
\hskip -2em
P^r_s=\delta^r_s-\frac{p_s\,\tilde\nabla^r\!H}{\Omega}=
\delta^r_s-\frac{v^r\,\tilde\nabla_{\!s}L}{\Omega}
\tag12.3
$$
(compare with formulas \thetag{10.12} and \thetag{1.5} above). Now
we see that \thetag{12.3} are components of operator-valued extended
tensor field $\bold P$ either in $\bold p$ and in $\bold
v$-representations respectively. It is easy to check up the
following identity:
$$
\hskip -2em
\bold P^2=\bold P\compos\bold P=\bold P.
\tag12.4
$$
Formula \thetag{12.4} means that $\bold P$ is a projector-valued
extended operator field in $M$. It projects along velocity vector
$\bold v$ onto a hyperplane defined by momentum covector $\bold p$.
When restricted to hypersurface $S$ due to lift \thetag{11.8} in
$\bold p$-representation or due to lift \thetag{11.9} in
$\bold v$-representation, it becomes projector field that projects
onto tangent hyperplane to $S$ along velocity vector $\bold v$.
In this restricted form components of $P$ appear in formula
\thetag{12.1}. This formula now can be written as
$$
\hskip -2em
\sum^n_{r=1}P^r_s\,\frac{\partial n_r}{\partial y^i}
=\frac{1}{\nu}\,\frac{\partial p_s}{\partial y^i}+\sum^n_{r=1}
\frac{1}{\nu}\left(\frac{1}{\Omega}\frac{\partial H}{\partial x^r}
-Q_r\!\right)\tau^r_i\,p_s.
\tag12.5
$$
Let's use the following identity that can be checked up immediately:
$$
\hskip -2em
\sum^n_{s=1}P^s_k\,p_s=\sum^n_{s=1}\left(\delta^s_k
-\frac{v^s\,p_k}{\Omega}\right)p_s=0.
\tag12.6
$$
Applying \thetag{12.4} and \thetag{12.6} to the equality \thetag{12.5},
we find that
$$
\hskip -2em
\sum^n_{r=1}P^r_s\,\frac{\partial n_r}{\partial y^i}=
\frac{1}{\nu}\,\sum^n_{r=1}P^r_s\,\frac{\partial p_r}{\partial y^i}.
\tag12.7
$$
Looking at \thetag{12.5}, we see that we cannot express partial
derivative $\partial n_s/\partial y^i$ in a pure form. However,
as appears, formula \thetag{12.7} is sufficient for our purposes.
Indeed, the term containing partial derivative $\partial n_r/
\partial y^j$ in \thetag{9.1} can be written as:
$$
\ldots
-\sum^n_{s=1}\sum^n_{r=1}\sum^n_{q=1}\sum^n_{k=1}\left(
\frac{p_q}{\Omega}\,\frac{\partial^2H}{\partial p_q\,
\partial p_k}\,Q_s+
\delta^k_q\,\frac{\partial Q_s}{\partial p_q}
\right)\nu^2\,P^r_k\frac{\partial n_r}{\partial y^j}
\,\tau^s_i+\ldots
$$
Therefore we can use the above identity \thetag{12.7} in order to
transform this term:
$$
\ldots
-\sum^n_{s=1}\sum^n_{r=1}\sum^n_{q=1}\sum^n_{k=1}\left(
\frac{p_q}{\Omega}\,\frac{\partial^2H}{\partial p_q\,
\partial p_k}\,Q_s+
\delta^k_q\,\frac{\partial Q_s}{\partial p_q}
\right)\nu\,P^r_k\frac{\partial p_r}{\partial y^j}
\,\tau^s_i+\ldots
$$
As a result of this transformation for the whole expression \thetag{9.1}
we obtain
$$
\hskip -0.1em
\gathered
\frac{\partial^2\nu}{\partial y^i\,\partial y^j}=\dots
+\sum^n_{s=1}\sum^n_{r=1}\left(\,\shave{\sum^n_{q=1}}
\shave{\sum^n_{k=1}}\frac{\nu}{\Omega^2}\,p_q\,p_k\,
\frac{\partial^2H}{\partial p_q\,\partial p_k}\,
\frac{\partial H}{\partial p_r}\,Q_s\,-\right.\\
\vspace{1ex}
\left.-\shave{\sum^n_{q=1}}\frac{\nu}{\Omega}\,p_q\,\frac{\partial^2H}
{\partial p_q\,\partial p_r}\,Q_s-\nu\,\frac{\partial Q_s}
{\partial p_r}+\shave{\sum^n_{q=1}}\frac{\nu}{\Omega}\,p_q\,
\frac{\partial H}{\partial p_r}\,\frac{\partial Q_s}{\partial p_q}
\right)\frac{\partial p_r}{\partial y^j}\,\tau^s_i\,+\\
\vspace{1ex}
+\sum^n_{s=1}\sum^n_{r=1}\left(-\shave{\sum^n_{q=1}}
\frac{\nu\,p_q}{\Omega}\,\frac{\partial Q_r}{\partial p_q}
\,\frac{\partial H}{\partial x^s}+\shave{\sum^n_{q=1}}
\frac{\nu\,p_q}{\Omega^2}\,\frac{\partial H}{\partial p_q}
\,\frac{\partial H}{\partial x^r}\,Q_s
+\shave{\sum^n_{q=1}}\shave{\sum^n_{k=1}}\frac{\nu\,p_q\,p_k}{\Omega^2}\,
\times\right.\\
\vspace{1ex}
\left.\times\,\frac{\partial^2H}{\partial p_q\,\partial p_k}\,
\frac{\partial H}{\partial x^r}\,Q_s
-\shave{\sum^n_{q=1}}\frac{\nu\,p_q}{\Omega}\,\frac{\partial^2H}
{\partial p_q\,\partial x^r}\,Q_s
+\shave{\sum^n_{q=1}}\nu\,p_q\,
\frac{\partial Q_r}{\partial p_q}\,Q_s
+\nu\,\frac{\partial Q_r}{\partial x^s}
\right)\tau^s_i\,\tau^r_j.
\endgathered\hskip -1em
\tag12.8
$$
Now we can replace partial derivative $\partial p_r/\partial y^j$
in \thetag{12.8} by covariant derivative $\nabla_{\!\boldsymbol
\tau_j}p_r$. For this purpose we shall use formula \thetag{11.2}
with $t=y^j$. This yields
$$
\hskip -2em
\frac{\partial p_r}{\partial y^j}=\nabla_{\!\boldsymbol\tau_j}p_r+
\sum^n_{\alpha=1}\sum^n_{\sigma=1}\Gamma^\alpha_{\sigma r}\,p_\alpha\,
\tau^\sigma_j.
\tag12.9
$$
Let's substitute \thetag{12.9} into \thetag{12.8} and let's use formulas
\thetag{10.3} and \thetag{10.7} in order to express partial derivatives
of $H$ and $Q_s$ through covariant derivatives:
$$
\gathered
\frac{\partial^2\nu}{\partial y^i\,\partial y^j}=\ldots
-\sum^n_{s=1}\sum^n_{r=1}\sum^n_{q=1}\left(p^q\,\frac{Q_s}{\Omega}
+\tilde\nabla^qQ_s\!\right)\nu\,P^r_q\,\nabla_{\!\boldsymbol
\tau_j}p_r\,\tau^s_i\,+\\
\vspace{1ex}
+\sum^n_{s=1}\sum^n_{r=1}\left(\frac{\nu\,|\bold p|^2\,\,
\nabla_{\!r}H}{\Omega^2}\,Q_s-\shave{\sum^n_{q=1}}p_q\,
\frac{\nu\,\,\nabla_{\!r}\!\tilde\nabla^qH}{\Omega}\,Q_s
-\nabla_rQ_s\,+\right.\\
\vspace{1ex}
+\left.
\shave{\sum^n_{q=1}}p_q\,\frac{\nu\,\,\nabla_{\!r}H}{\Omega}\,
\tilde\nabla^qQ_s+\frac{\nu\,\,\nabla_{\!r}H}{\Omega}\,Q_s+
\shave{\sum^n_{q=1}}\nu\,p_q\,Q_s\,\tilde\nabla^qQ_r
\right)\tau^s_i\,\tau^r_j.
\endgathered
\hskip -1em
\tag12.10
$$
By transforming \thetag{9.2} we can obtain another expression for
the same derivative:
$$
\gathered
\frac{\partial^2\nu}{\partial y^j\,\partial y^i}=\ldots
-\sum^n_{s=1}\sum^n_{r=1}\sum^n_{q=1}\left(p^q\,\frac{Q_r}{\Omega}
+\tilde\nabla^qQ_r\!\right)\nu\,P^s_q\,\nabla_{\!\boldsymbol
\tau_i}p_s\,\tau^r_j\,+\\
\vspace{1ex}
+\sum^n_{s=1}\sum^n_{r=1}\left(\frac{\nu\,|\bold p|^2\,\,
\nabla_{\!s}H}{\Omega^2}\,Q_r-\shave{\sum^n_{q=1}}p_q\,
\frac{\nu\,\,\nabla_{\!s}\!\tilde\nabla^qH}{\Omega}\,Q_r
-\nabla_sQ_r\,+\right.\\
\vspace{1ex}
+\left.
\shave{\sum^n_{q=1}}p_q\,\frac{\nu\,\,\nabla_{\!s}H}{\Omega}\,
\tilde\nabla^qQ_r+\frac{\nu\,\,\nabla_{\!s}H}{\Omega}\,Q_r+
\shave{\sum^n_{q=1}}\nu\,p_q\,Q_r\,\tilde\nabla^qQ_s
\right)\tau^s_i\,\tau^r_j.
\endgathered
\hskip -1em
\tag12.11
$$
Here, as in section~9, we denote by dots terms symmetric in
indices $r$ and $s$. They do not affect ultimate compatibility
equation, which will be derived from \thetag{12.10} and
\thetag{12.11}. In the equalities \thetag{12.10} and \thetag{12.11}
we used matrix $g^{ij}=\tilde\nabla^i\tilde\nabla^jH$ from
\thetag{10.13} as metric tensor and introduced the following
notations:
$$
\xalignat 2
&\hskip -2em
p^q=\sum^n_{k=1}g^{qk}\,p_q,
&&|\bold p|^2=\sum^n_{q=1}\sum^n_{k=1}g^{qk}\,p_q\,p_k.
\tag12.12
\endxalignat
$$
In general, matrix \thetag{10.13} is not positive, therefore
$|\bold p|^2$ in \thetag{12.12} is not necessarily positive
number for $\bold p\neq 0$.\par
     Let's equate right hand sides of \thetag{12.10} and
\thetag{12.11}. As a result we obtain compatibility equation
\thetag{7.19} written in terms of covariant derivatives:
$$
\gathered
\aligned
\sum^n_{s=1}\sum^n_{r=1}\sum^n_{q=1}&\left(p^q\,\frac{Q_s}{\Omega}
+\tilde\nabla^qQ_s\!\right)P^r_q\,\nabla_{\!\boldsymbol
\tau_j}p_r\,\tau^s_i\,-\\
\vspace{1ex}
&-\sum^n_{s=1}\sum^n_{r=1}\sum^n_{q=1}\left(p^q\,\frac{Q_r}{\Omega}
+\tilde\nabla^qQ_r\!\right)P^s_q\,\nabla_{\!\boldsymbol
\tau_i}p_s\,\tau^r_j=\\
\endaligned
\vspace{1ex}
=\sum^n_{s=1}\sum^n_{r=1}\left(\frac{|\bold p|^2\,\,
\nabla_{\!r}H}{\Omega^2}\,Q_s-\shave{\sum^n_{q=1}}p_q\,
\frac{\nabla_{\!r}\!\tilde\nabla^qH}{\Omega}\,Q_s
-\nabla_rQ_s\,+\right.\\
\vspace{1ex}
+\left.
\shave{\sum^n_{q=1}}p_q\,\frac{\nabla_{\!r}H}{\Omega}\,
\tilde\nabla^qQ_s+\frac{\nabla_{\!r}H}{\Omega}\,Q_s+
\shave{\sum^n_{q=1}}p_q\,Q_s\,\tilde\nabla^qQ_r
\right)\tau^s_i\,\tau^r_j-\\
\vspace{1ex}
-\sum^n_{s=1}\sum^n_{r=1}\left(\frac{|\bold p|^2\,\,
\nabla_{\!s}H}{\Omega^2}\,Q_r-\shave{\sum^n_{q=1}}p_q\,
\frac{\nabla_{\!s}\!\tilde\nabla^qH}{\Omega}\,Q_r
-\nabla_sQ_r\,+\right.\\
\vspace{1ex}
+\left.
\shave{\sum^n_{q=1}}p_q\,\frac{\nabla_{\!s}H}{\Omega}\,
\tilde\nabla^qQ_r+\frac{\nabla_{\!s}H}{\Omega}\,Q_r+
\shave{\sum^n_{q=1}}p_q\,Q_r\,\tilde\nabla^qQ_s
\right)\tau^s_i\,\tau^r_j.
\endgathered
\tag12.13
$$
The most important point here is that the form of equation
\thetag{12.13} does not depend on which particular connection
is used in covariant derivatives. One can choose any symmetric
extended connection $\Gamma$ in $M$.
\head
13. Geometry of hypersurfaces.
\endhead
    Though most terms in compatibility equation \thetag{12.13}
can be treated as components of extended tensor fields in $M$,
the equation in whole is related to hypersurface $S$. This
reveals if we look at partial derivatives $\nabla_{\!\boldsymbol
\tau_j}p_r$ and $\nabla_{\!\boldsymbol\tau_i}p_s$. For further
treatment of \thetag{12.13} we should learn how to treat these
partial derivatives.\par
   Let $\boldsymbol\tau$ be some arbitrary tangent vector to $S$.
It represented as linear combination of basic tangent vectors:
$\boldsymbol\tau=\alpha^1\cdot\boldsymbol\tau_1+\ldots+
\alpha^{n-1}\cdot\boldsymbol\tau_{n-1}$. Let
$$
\hskip -2em
f(\boldsymbol\tau)=\nabla_{\!\boldsymbol\tau}\bold p=
\sum^n_{r=1}\sum^{n-1}_{j=1}\left(\alpha^j\,
\nabla_{\!\boldsymbol\tau_j}p_r\right)
\kern -1pt\cdot dx^r.
\tag13.1
$$
Then \thetag{13.1} defines linear map $f\!:T_p(S)\to T^*_p(M)$.
Let's consider composite map
$$
\hskip -2em
\bold b=-\bold P^*\compos f\compos\bold P.
\tag13.2
$$
Projection operator $\bold P^*$ in \thetag{13.2} is a conjugate
operator for projector $\bold P$ with components \thetag{12.2}.
Remember that $\bold P$ projects onto the subspace $T_p(S)$ in
$T_p(M)$ for $p\in S$. Therefore linear map $\bold b\!:T_p(M)\to
T^*_p(M)$ is correctly defined by formula \thetag{13.2}. Each map
from $T_p(M)$ to conjugate space $T^*_p(S)$ is given by some
bilinear form. In our case this bilinear form is defined by formula
$$
\hskip -2em
b(\bold X,\bold Y)=\left<\bold b(\bold Y)\,|\,\bold X\right>.
\tag13.3
$$
Due to the presence of projection operators $\bold P$ and $\bold P^*$
in \thetag{13.2} we have
$$
\hskip -2em
b(\bold X,\bold Y)=b(\bold P(\bold X),\bold Y)=b(\bold X,
\bold P(\bold Y)).
\tag13.4
$$
\proclaim{Theorem 13.1} Bilinear form \thetag{13.3} defined by
linear map \thetag{13.2} is symmetric.
\endproclaim
   Bilinear form \thetag{13.3} restricted to tangent space $T_p(S)$
of hypersurface $S$ is known as {\it second fundamental form} of
hypersurface $S$. Its components
$$
\hskip -2em
\beta_{ij}=b(\boldsymbol\tau_i,\boldsymbol\tau_j)
\tag13.5
$$
define tensor field in inner geometry of hypersurface $S$. But
for our purposes coordinate representation of bilinear form
\thetag{13.3} in outer geometry is more preferable:
$$
b=\sum^n_{i=1}\sum^n_{j=1}b_{ij}\,dx^i\,dx^j.
$$
In order to prove theorem~13.1 we need some preliminary results.
First is given by the following lemma for extended connection
components.
\proclaim{Lemma 13.1} For any symmetric extended connection
$\Gamma$ in $M$ and for any fixed point $q_0=(p_0,\bold p)$
of cotangent bundle $T^*\!M$ there is local chart of $M$ such
that all connection components $\Gamma^k_{ij}(q)=\Gamma^k_{ij}
(x^1,\ldots,x^n,p_1,\ldots,p_n)$ in this chart become zero
at the point $q=q_0$.
\endproclaim
Lemma~13.1 is formulated for extended connection in
$\bold p$-representation and for its components \thetag{10.6}.
Similar lemma can be stated for extended connection in
$\bold v$-representation and for its components
\thetag{10.15}.
\proclaim{Lemma~13.2} For any symmetric extended connection
$\Gamma$ in $M$ and for any fixed point $q_0=(p_0,\bold v)$
of tangent bundle $TM$ there is local chart of $M$ such
that all connection components $\Gamma^k_{ij}(q)=\Gamma^k_{ij}
(x^1,\ldots,x^n,v^1,\ldots,v^n)$ in this chart become zero
at the point $q=q_0$.
\endproclaim
    Here one should emphasize that lemmas~13.1 and 13.2 assert
vanishing of $\Gamma^k_{ij}$ only at single point $q=q_0$. Sometimes
they can vanish in whole neighborhood of this point, but this is
not the case of general position. Proof of lemmas~13.1 and 13.2
is rather standard. It is based upon formula \thetag{10.4}. One
can find proof of analogous lemma in Chapter~\uppercase
\expandafter{\romannumeral 5} of thesis \cite{6}. With minor changes
this proof is applicable to lemmas~13.1 and 13.2.\par
    Now suppose that $x^1,\,\ldots,\,x^n$ are coordinates in a chart,
existence of which is asserted by lemma~13.1.
Then, using transformation formula \thetag{10.4}, one can check up
that linear change variables (i\.\,e\. that, for which transition
matrices $S$ and $T$ in \thetag{10.5} are constant) preserves the
property of vanishing of $\Gamma^k_{ij}(q)$ at the point $q_0=(p_0,
\bold p)$. If we apply this fact to a fixed point $p_0\in S$, we
can choose local coordinates $x^1,\,\ldots,\,x^n$ in some neighborhood
of $p_0$ such that
\roster
\rosteritemwd=1pt
\item"1)" fixed point $p_0$ is the origin in coordinates
$x^1,\,\ldots,\,x^n$,
i\.\,e\. $x^i(p_0)=0$ for the whole range of index $i=1,\,\ldots,
\,n$;
\item"2)" hypersurface $S$ is given parametrically by the equations
$$
\hskip -1em
\cases
x^1=y^n,\\
\vspace{-1ex}
.\ .\ .\ .\ .\ .\ .\ .\ \\
x^{n-1}=y^{n-1},\\
\vspace{0.5ex}
x^n=z(y^1,\ldots,y^{n-1}),
\endcases
\tag13.6
$$
where $z(y^1,\ldots,y^{n-1})$ is some smooth function of parameters
$y^1,\,\ldots,\,y^{n-1}$ vanishing at the origin, i\.\,e\.
$z(0,\ldots,0)=0$, and having extremum there;
\item"3)" normal covector $\bold n=\bold n(y^1,\ldots,y^{n-1})$ of
hypersurface $S$ related with momentum covector $\bold p$ on $S$ by
the equality \thetag{11.8} is given by its components:
$$
\hskip -1em
\bold n=(-z'_1,\ \ldots, -z'_n, 1)\text{, \ where \ }
z'_i=\frac{\partial z}{\partial y^i};
\tag13.7
$$
\item"4)" connection components $\Gamma^k_{ij}=\Gamma^k_{ij}(p,
\bold p(p))$, where $\bold p=\bold p(p)$ is determined by the
equality \thetag{11.8}, do vanish at the origin $p=p_0$.
\endroster
Formula \thetag{13.7} is derived directly from \thetag{13.6}.
Indeed, if we calculate tangent vectors $\boldsymbol\tau_1,\,
\ldots,\,\boldsymbol\tau_{n-1}$ by using formula \thetag{2.5},
for $\tau^s_i$ we obtain
$$
\hskip -2em
\tau^s_i=\cases 0 &\text{for \ }s\neq i, n,\\
1 &\text{for \ }s=i,\\
z'_i &\text{for \ }i=n.
\endcases
\tag13.8
$$
Now it's easy to see that covector \thetag{13.7} is orthogonal to
vectors $\boldsymbol\tau_1,\,\ldots,\,\boldsymbol\tau_{n-1}$ in
the sense of the equality \thetag{7.3}. As we know, normal covector
of hypersurface is determined up to a scalar factor. This uncertainty
in \thetag{13.6} is eliminated by the condition that last component
of $\bold n$ is equal to unity.\par
   Let's consider components of covector $\nabla_{\!\boldsymbol\tau}
\bold p$ in \thetag{13.1} for those special coordinates $x^1,\,\ldots,
\,x^n$ we have chosen above. Applying \thetag{11.2} and \thetag{12.7},
we derive
$$
\hskip -2em
\sum^n_{r=1}P^r_s\,\nabla_{\!\boldsymbol\tau_i}p_r=
\sum^n_{r=1}\nu\,P^r_s\,\nabla_{\!\boldsymbol\tau_i}n_r.
\tag13.9
$$
Applying formula \thetag{11.2} again, for covariant derivative
$\nabla_{\!i}n_r$ we obtain
$$
\hskip -2em
\nabla_{\!\boldsymbol\tau_i}n_r=\frac{\partial n_r}{\partial y^i}
-\sum^n_{\alpha=1}\sum^n_{\sigma=1}\Gamma^\alpha_{\sigma r}
\,n_\alpha\,\tau^\sigma_j.
\tag13.10
$$
Remember that in the above formula $\Gamma^\alpha_{\sigma r}=0$
for $p=p_0$. Therefore, taking into account formula \thetag{13.7},
from \thetag{13.10} we derive the following equality:
$$
\hskip -2em
\nabla_{\!\boldsymbol\tau_i}n_r\,\hbox{\vrule height 8pt depth
8pt width 0.5pt}_{\,p=p_0}=\cases z''_{ir} &\text{for \ }r<n,\\
0 &\text{for \ }r=n.\endcases
\tag13.11
$$
Projector $\bold P$ projects onto the hyperplane $T_p(S)$ and
$\boldsymbol\tau_i\in T_p(S)$, hence $\bold P(\boldsymbol\tau_i)
=\boldsymbol\tau_i$. Then for $\bold b(\boldsymbol\tau_i)$,
applying \thetag{13.1}, \thetag{13.2}, and \thetag{13.9}, we
derive
$$
\hskip -2em
\bold b(\boldsymbol\tau_i)\,\hbox{\vrule height 8pt depth 8pt
width 0.5pt}_{\,p=p_0}=-\sum^n_{s=1}\sum^n_{r=1}\left(\nu\,P^r_s\
\nabla_{\!\boldsymbol\tau_i}n_r\right)\kern -1pt\cdot dx^s.
\tag13.12
$$
Further, using formulas \thetag{13.3}, \thetag{13.5}, \thetag{13.11},
and \thetag{13.12}, we obtain
$$
\hskip -2em
\beta_{ij}\,\hbox{\vrule height 8pt depth 8pt
width 0.5pt}_{\,p=p_0}=-\sum^n_{s=1}\sum^n_{r=1}\nu\,P^r_s\
\nabla_{\!\boldsymbol\tau_i}n_r\,\tau^s_j=-\sum^n_{s=1}
\sum^n_{r=1}\nu\,P^r_s\,z''_{ir}\,\tau^s_j.
\tag13.13
$$
Now remember formula \thetag{12.2} for $P^r_s$ and formula
\thetag{13.7} for $\bold n$. Moreover, let's remember that
function $z(y^1,\ldots,y^{n-1})$ has extremum at the origin.
This means that $z'_i=0$ for $p=p_0$. Then from \thetag{12.2},
\thetag{13.7}, and \thetag{13.8} we get
$$
\xalignat 2
&\hskip -2em
P^r_s\,\hbox{\vrule height 8pt depth 8pt width 0.5pt}_{\,p=p_0}
=\delta^r_s-\frac{v^r\,\delta^n_s}{v^n},
&&\tau^s_i\,\hbox{\vrule height 8pt depth 8pt width 0.5pt}_{\,p=p_0}
=\delta^s_i.
\tag13.14
\endxalignat
$$
Substituting \thetag{13.14} into \thetag{13.13} we now obtain
the following equality:
$$
\hskip -2em
\beta_{ij}\,\hbox{\vrule height 8pt depth 8pt
width 0.5pt}_{\,p=p_0}=-\nu\,z''_{ij}
=-\nu\,\frac{\partial^2 z}{\partial y^i
\,\partial y^j}.
\tag13.15
$$
Looking at \thetag{13.15}, it's easy to see that second fundamental
form of hypersurface $S$ is symmetric: $\beta_{ij}=\beta_{ij}$. This
is in concordance with classical results for hypersurfaces in
Riemannian manifolds. But now we are in quite different geometry.
\par
    Second fundamental form $\boldsymbol\beta$ with components given
by formula \thetag{13.5} is a tensor field in $S$. It is obtained by
restricting quadratic form \thetag{13.3} to inner geometry of $S$.
Let's calculate components of this form in outer geometry, i\.\,e\.
in local chart of manifold $M$ at the point $p=p_0$. Due to \thetag{13.14}
tangent vectors $\boldsymbol\tau_1,\,\ldots,\,\boldsymbol\tau_{n-1}$
coincide with $n-1$ coordinate tangent vectors at the point $p_0$:
$$
\hskip -2em
\boldsymbol\tau_i\,\hbox{\vrule height 8pt depth 8pt width
0.5pt}_{\,p=p_0}=\bold E_i=\frac{\partial}{\partial x^i}
\text{\ \ for \ }i=1,\,\ldots,\,n-1.
\tag13.16
$$
As a consequence of \thetag{13.16}, \thetag{13.3}, and \thetag{13.5}
we obtain
$$
\hskip -2em
b_{ij}=\bold b(\bold E_i,\bold E_j)=\bold b(\boldsymbol\tau_i,
\boldsymbol\tau_j)=\beta_{ij}\text{\ \ for \ }1\leqslant i,j<n.
\tag13.17
$$
Let's take $n$-th coordinate vector $\bold E_n=\partial/\partial x^n$.
From \thetag{13.14} we derive
$$
\bold P(\bold E_n)\,\hbox{\vrule height 8pt depth 8pt width
0.5pt}_{\,p=p_0}=\sum^n_{i=1}P^r_n\cdot\bold E_r=-\sum^{n-1}_{r=1}
\frac{v^i}{v_n}\cdot\bold E_r=-\sum^{n-1}_{i=1}\frac{v^i}{v_n}\cdot
\boldsymbol\tau_r.
$$
Using the above equality and formulas \thetag{13.4} and \thetag{13.5},
for $i<n$ we find
$$
\hskip -2em
b_{in}=b(\bold E_i,\bold E_n)=b(\bold E_i,\bold P(\bold E_n))=
-\sum^{n-1}_{r=1}\frac{\beta_{ir}\,v^r}{v^n}.
\tag13.18
$$
In a similar way for component $b_{ni}$ of quadratic form
\thetag{13.3} we obtain
$$
\hskip -2em
b_{ni}=b(\bold E_n,\bold E_i)=b(\bold P(\bold E_n),\bold E_i)=
-\sum^{n-1}_{r=1}\frac{\beta_{ri}\,v^r}{v^n}.
\tag13.19
$$
Comparing \thetag{13.18} with \thetag{13.19} and taking into account
symmetry of $\beta_{ij}$ (see \thetag{13.15} above), we find that
$b_{in}=b_{ni}$ for $i<n$. Combining this equality with \thetag{13.17}
and taking into account symmetry of $\beta_{ij}$ again, we get
$$
\hskip -2em
b_{ij}=b_{ji}
\tag13.20
$$
for all $i$ and $j$. We have proved the equality \thetag{13.20} at one
point $p=p_0$ by choosing special local chart in $M$. But we can choose
arbitrary point of $M$ for $p_0$. Besides, $b_{ij}$ are components of
tensor. Therefore, being symmetric in one local chart, they keep symmetry
when transformed to another chart. This proves theorem~13.1.
\proclaim{Theorem 13.2} Let $q_0=(p_0,\bold p)$ be some fixed point
of cotangent bundle $T^*\!M$ with $\bold p\neq 0$ and let projector
$\bold P$ be the value of projector-valued extended tensor field
\thetag{12.4} at this point. Then any symmetric quadratic form $b$ in
$T_{p_0}(M)$ satisfying the equality \thetag{13.4} can be determined by
some hypersurface $S$ passing through the point $p_0$ and tangent to
null-space of covector $\bold p$ at this point.
\endproclaim
    Proof of theorem~13.2 now is very simple. Indeed, quadratic form
$b$ satisfying the equality \thetag{13.4} is completely determined by
its restriction to null-space of covector $\bold p$. Then we can choose
local chart in $M$ with $\Gamma^k_{ij}(p_0,\bold p)=0$ and with first
$n-1$ coordinate vectors $\bold E_1,\,\ldots,\,\bold E_{n-1}$ all being
in null-space of covector $\bold p$ at the point $p_0$. Let's define
matrix $\beta$ by means of components of quadratic form $b$:
$$
\hskip -2em
\beta_{ij}=b_{ij}=b(\bold E_i,\bold E_j)\,\hbox{\vrule height 8pt
depth 8pt width 0.5pt}_{\,p=p_0}\text{\ \ for \ }
1\leqslant i,j<n.
\tag13.21
$$
Now it's sufficient to define hypersurface $S$ by parametric equations
\thetag{13.6} and take the following function $z=z(y^1,\,\ldots,y^{n-1})$
in them:
$$
\hskip -2em
z=-\frac{1}{2\,\nu_0}\sum^n_{i=1}\sum^n_{j=1}
\beta_{ij}\,y^i\,y^j.
\tag13.22
$$
Here $\nu_0=\nu(p_0)$ is a constant from \thetag{7.15}. Due to
equalities \thetag{13.15}, \thetag{13.17}, and \thetag{13.21}
it is clear that such hypersurface $S$ reproduces quadratic form $b$
used to define it through \thetag{13.21} and \thetag{13.22}. Thus,
theorem~13.2 is proved.\par
\head
14. Additional normality equations.
\endhead
    Now let's apply theorem~13.2 to the study of the equation
\thetag{12.13}. Let's fix some point $p_0$ in $M$ and some covector
$\bold p$ at this point. Then, relying upon theorem~13.2, let's take
hypersurface $S$ passing through this point tangent to null-space
of covector $\bold p$ and such that its second fundamental form
$\boldsymbol\beta$ is zero at the point $p=p_0$. This means that
quadratic form $b$ in \thetag{13.3} is also zero for $p=p_0$. Then
due to \thetag{13.1} and \thetag{13.2} we get $\nabla_{\!\boldsymbol
\tau_j}p_r=0$. Equivalently, $\nabla_{\!\boldsymbol\tau_i}p_s=0$.
Substituting these two equalities into \thetag{12.13}, we find that
first two sums in \thetag{12.13} do vanish. Looking at other terms,
we see that they are components of an extended tensor field of type
$(0,2)$ contracted with two vectors $\boldsymbol\tau_i$ and
$\boldsymbol\tau_j$ tangent to $S$. These two vectors depend on the
choice of parameters $y^1,\,\ldots,\,y^n$, i\.\,e\. on the choice
of local chart on $S$. By choosing this local chart properly we can
associate $\boldsymbol\tau_i$ and $\boldsymbol\tau_j$ with two
arbitrary vectors in $T_{p_0}(S)$. This means that we can write
$$
\xalignat 2
&\hskip -2em
\boldsymbol\tau_i=\bold P(\bold X),
&&\boldsymbol\tau_j=\bold P(\bold Y),
\tag14.1
\endxalignat
$$
where $\bold X$ and $\bold Y$ are two arbitrary vectors in $T_{p_0}
(M)$. Due to arbitrariness of vectors $\bold X$ and $\bold Y$
in \thetag{14.1} from rest part of \thetag{12.13} we derive
$$
\hskip -2em
\gathered
\sum^n_{s=1}\sum^n_{r=1}\left(\frac{|\bold p|^2\,\,
\nabla_{\!r}H}{\Omega^2}\,Q_s-\shave{\sum^n_{q=1}}p_q\,
\frac{\nabla_{\!r}\!\tilde\nabla^qH}{\Omega}\,Q_s
-\nabla_rQ_s\,+\right.\\
\vspace{1ex}
+\left.
\shave{\sum^n_{q=1}}p_q\,\frac{\nabla_{\!r}H}{\Omega}\,
\tilde\nabla^qQ_s+\frac{\nabla_{\!r}H}{\Omega}\,Q_s+
\shave{\sum^n_{q=1}}p_q\,Q_s\,\tilde\nabla^qQ_r
\right)P^s_i\,P^r_j=\\
\vspace{1ex}
=\sum^n_{s=1}\sum^n_{r=1}\left(\frac{|\bold p|^2\,\,
\nabla_{\!s}H}{\Omega^2}\,Q_r-\shave{\sum^n_{q=1}}p_q\,
\frac{\nabla_{\!s}\!\tilde\nabla^qH}{\Omega}\,Q_r
-\nabla_sQ_r\,+\right.\\
\vspace{1ex}
+\left.
\shave{\sum^n_{q=1}}p_q\,\frac{\nabla_{\!s}H}{\Omega}\,
\tilde\nabla^qQ_r+\frac{\nabla_{\!s}H}{\Omega}\,Q_r+
\shave{\sum^n_{q=1}}p_q\,Q_r\,\tilde\nabla^qQ_s
\right)P^s_i\,P^r_j.
\endgathered
\tag14.2
$$
Note that the equations \thetag{14.2} are partial differential
equations for components of extended covector field $\bold Q$,
which is used in the equations of Newtonian dynamics, when they
are written in the form \thetag{7.5} relative to modified
Hamiltonian dynamical. They are written in terms of covariant
derivatives \thetag{10.3} and \thetag{10.7}. Though we used some
special hypersurface $S$ in order to derive them, in their ultimate
form they do not depend on any particular choice of $S$ and,
moreover, they can be written in the absence of $S$ at all.\par
    Now remember that $\bold P(\boldsymbol\tau_i)=\boldsymbol\tau_i$
and $\bold P(\boldsymbol\tau_j)=\boldsymbol\tau_j$. Therefore we
can apply \thetag{14.2} back to \thetag{12.13}. As a result we obtain
the following equality:
$$
\hskip -2em
\aligned
\sum^n_{s=1}\sum^n_{r=1}\sum^n_{q=1}&\left(p^q\,\frac{Q_s}{\Omega}
+\tilde\nabla^qQ_s\!\right)P^r_q\,\nabla_{\!\boldsymbol
\tau_j}p_r\,\tau^s_i\,=\\
\vspace{1ex}
&=\sum^n_{s=1}\sum^n_{r=1}\sum^n_{q=1}\left(p^q\,\frac{Q_r}{\Omega}
+\tilde\nabla^qQ_r\!\right)P^s_q\,\nabla_{\!\boldsymbol
\tau_i}p_s\,\tau^r_j.
\endaligned
\tag14.3
$$
This means that the equation \thetag{12.13} splits into two parts,
first part is \thetag{14.3}, while second part leads to the equations
\thetag{14.2}.\par
     Both sides of \thetag{14.3} do vanish for our special hypersurface
$S$ above at its fixed point $p=p_0$. However, for arbitrary hypersurface
$S$ they are nonzero, therefore we are to study \thetag{14.3} in order to
derive other equations for covector field $\bold Q$. Let's denote by
$\bold B$ extended tensor field with components
$$
\hskip -2em
B^r_s=\sum^n_{k=1}\sum^n_{q=1}P^r_q\left(p^q\,\frac{Q_k}
{\Omega}+\tilde\nabla^qQ_k\!\right)P^k_s.
\tag14.4
$$
It's clear that tensor field $\bold B$ with components \thetag{14.4}
is an operator field. Due to the presence of $P^k_s$ and $P^r_q$ in
\thetag{14.4} we have the following equalities:
$$
\hskip -2em
\bold B=\bold P\compos\bold B=\bold B\compos\bold P.
\tag14.5
$$
Relying upon \thetag{13.1}, \thetag{13.2}, \thetag{13.3} and
using \thetag{14.4}, now we can write \thetag{14.3} as
$$
\hskip -2em
b(\boldsymbol\tau_i,\bold B\boldsymbol\tau_j)=
b(\boldsymbol\tau_j,\bold B\boldsymbol\tau_i).
\tag14.6
$$
As we noted above, vectors $\boldsymbol\tau_i$ and $\boldsymbol\tau_j$
can be replaced by two arbitrary vectors $\bold X$ and $\bold Y$, see
formulas \thetag{14.1}. Then \thetag{14.6} is transformed to
$$
b(\bold X,\bold B\compos\bold P(\bold Y))=
b(\bold Y,\bold B\compos\bold P(\bold X)).
$$
Due to \thetag{14.5} and theorem~13.1 we can further simplify this
relationship:
$$
\hskip -2em
b(\bold X,\bold B(\bold Y))=b(\bold B(\bold X),\bold Y).
\tag14.7
$$
Formula \thetag{14.7} means that operator $\bold B$ is symmetric with
respect to bilinear form \thetag{13.3}. Now we are to utilize this
equality.\par
    Let's denote by $W=\Img\bold P$ the image of projection operator
$\bold P=\bold P(q_0)$ for some fixed point $q_0=(p_0,\bold p)$ of
cotangent bundle $T^*\!M$ (see theorem~13.2). Then $W$ is
$(n-1)$-dimensional subspace in $T_{p_0}(M)$. It coincides with
null-space of covector $\bold p$. Due to \thetag{14.5} subspace $W$
is invariant under the action of operator $\bold B$. Moreover,
operator $\bold B$ is completely determined by its restriction to $W$.
For instance, if restriction of $\bold B$ to $W$ is identical operator
in $W$, then $\bold B=\bold P$:
$$
\hskip -2em
\bold B\,\hbox{\vrule height 8pt depth 6pt width 0.5pt}_{\,W}=
\id_W\text{\ \ implies \ }\bold B=\bold P.
\tag14.8
$$\par
Bilinear form $b$ is also completely determined by its restriction to
subspace $W$. This follows from relationships \thetag{13.4}. The
equality \thetag{14.7} means that restriction of operator $B$ to $W$
is symmetric with respect to restriction of $b$ to $W$. Now recall
theorem~13.2. It means that arbitrary quadratic form in subspace $W$
can be obtained as second fundamental form of some hypersurface $S$
tangent to $W$. Hence restriction of $\bold B$ to $W$ is symmetric
with respect to all quadratic forms in $W$. It takes place if and
only if the restriction of $\bold B$ to $W$ is a scalar operator,
i\.\,e\.
$$
\bold B\,\hbox{\vrule height 8pt depth 6pt width 0.5pt}_{\,W}=
\lambda\cdot\id_W.
$$
This is easy result in linear algebra. Now, applying \thetag{14.8}
to the above equality, for operator $\bold B$ itself we derive the
following representation:
$$
\hskip -2em
\bold B=\lambda\cdot\bold P.
\tag14.9
$$
Here $\lambda$ is some scalar factor. It can be expressed explicitly
through trace of $\bold B$:
$$
\hskip -2em
\lambda=\frac{\tr\bold B}{n-1}.
\tag14.10
$$
Formula \thetag{14.9} is important result. Combining it with
\thetag{14.4}, we obtain
$$
\hskip -2em
\sum^n_{k=1}\sum^n_{q=1}P^r_q\left(p^q\,\frac{Q_k}{\Omega}
+\tilde\nabla^qQ_k\!\right)P^k_s=\lambda\,P^r_s.
\tag14.11
$$
Substituting \thetag{14.10} for scalar factor $\lambda$
in \thetag{14.11}, we get the following equality:
$$
\hskip -2em
\aligned
\sum^n_{k=1}\sum^n_{q=1}&P^r_q\left(p^q\,\frac{Q_k}{\Omega}
+\tilde\nabla^qQ_k\!\right)P^k_s=\\
\vspace{1ex}
&=\sum^n_{k=1}\sum^n_{q=1}\left(p^q\,\frac{Q_k}
{\Omega}+\tilde\nabla^qQ_k\!\right)\frac{P^k_q\,P^r_s}{n-1}.
\endaligned
\tag14.12
$$
This is another additional normality equation. Both \thetag{14.2} and
\thetag{14.12} form a system of partial differential equations for
components of extended convector field $\bold Q$ that defines Newtonian
dynamical system in form \thetag{7.5} relative to Hamiltonian dynamical
system with Hamilton function $H$. Having derived these equations, we
proved the following theorem.
\proclaim{Theorem 14.1} Additional normality condition stated in
definition~7.2 for Newtonian dynamical system \thetag{7.5} in
multidimensional case $n\geqslant 3$ is equivalent to the system of
additional normality equations \thetag{14.2} and \thetag{14.12} that
should be fulfilled at all points $q=(p,\bold p)$ of cotangent bundle
$T^*\!M$, where $\bold p\neq 0$.
\endproclaim
\head
15. Connection invariance.
\endhead
    Deriving additional normality equations \thetag{14.2} and
\thetag{14.12} above we used some symmetric extended connection
$\Gamma$. Components of this connection are present in \thetag{14.2}
due to covariant derivatives $\nabla_r$ and $\nabla_s$. However,
above we did not specify which particular symmetric connection $\Gamma$
is used. This means that additional normality equations \thetag{14.2}
and \thetag{14.12} should be invariant under transformations
$$
\hskip -2em
\Gamma^k_{ij}\to \Gamma^k_{ij}+T^k_{ij},
\tag15.1
$$
where $T^k_{ij}$ are components of symmetric extended tensor field
of type $(1,2)$. Differential equations \thetag{14.12} are obviously
invariant under transformations \thetag{15.1} since momentum gradient
$\tilde\nabla$ is defined without use of connection components (see
formula \thetag{10.3}). In general, differential equations \thetag{14.2}
are not invariant under these transformations. However, if component
of covector field $\bold Q$ satisfy differential equations \thetag{14.12},
then equations \thetag{14.2} become invariant under transformations
\thetag{15.1}.
In other words, this means that differential equations \thetag{14.2} are
invariant under transformations \thetag{15.1} modulo differential
equations \thetag{14.12}. This fact can be checked up by direct
calculations.
\head
16. Weak normality equations.
\endhead
    Let's fix some point $q_0=(p_0,\bold p_0)$ of cotangent bundle
$T^*\!M$ with $\bold p_0\neq 0$. It yields initial data for Newtonian
dynamical system written in form \thetag{7.5} and defines a trajectory
$p=p(t)$ of this dynamical system passing through the point $p=p_0$ at
time instant $t=0$. Null-space of covector $\bold p_0$ is a hyperplane
in tangent space $T_{p_0}(M)$. Let's denote it by $\alpha$. One can
draw various hypersurfaces passing through the point $p_0$ and tangent
to hyperplane $\alpha$ at this point. Suppose that $S$ is one of such
hypersurfaces and suppose that $\bold n=\bold n(p)$ is smooth normal
covector of $S$ in some neighborhood of the point $p_0$. At the very
point $p=p_0$ we have the equality
$$
\hskip -2em
\bold p_0=\nu_0\cdot\bold n(p_0)\text{, \ where \ }\nu_0\neq 0.
\tag16.1
$$
Now let's take some smooth function on $S$ normalized by the condition
\thetag{7.15} and set up initial data \thetag{7.4} for Newtonian dynamical
system \thetag{7.5}. Solving Cauchy problem with these initial data, we
get a family of trajectories for dynamical system \thetag{7.5}, which
includes our initial trajectory passing through the point $p=p_0$. This
is easily seen if we compare \thetag{16.1} and \thetag{7.15}.\par
    In local coordinates the family of trajectories constructed just
above is represented by functions \thetag{2.4}. Using them, we define
variation vectors \thetag{2.5} and deviation functions \thetag{2.6}.
Now suppose that Newtonian dynamical system \thetag{7.5} satisfy strong
normality condition (see definition~8.1). This means that by proper
choice of function $\nu(p)$ we can make all deviation functions
$\varphi_1,\,\ldots,\,\varphi_{n-1}$ to be identically zero. Hence
initial conditions \thetag{7.2} are fulfilled. From section~7 we
know that initial data \thetag{7.2} are equivalent to Pfaff equations
\thetag{7.13} for $\nu$. We can vary constant $\nu_0\neq 0$ in
normalizing condition \thetag{7.15}, and for each value of this constant
due to strong normality condition we would have some function $\nu(p)$
on $S$ satisfying Pfaff equations \thetag{7.13}. Due to lemma~7.2 then
Pfaff equations \thetag{7.13} are compatible in the sense of
definition~7.1. Thus we have proved the following theorem.
\proclaim{Theorem 16.1} Strong normality condition implies additional
normality condition for Newtonian dynamical system \thetag{7.5}.
\endproclaim
Additional normality condition is formulated in definition~7.2. In
section~14 we have shown that additional normality condition is
equivalent to additional normality equations \thetag{14.2} and
\thetag{14.12}. Thus strong normality condition leads to the equations
\thetag{14.2} and \thetag{14.12} for components of covector field
$\bold Q$. However, it can yield much more. Indeed, it implies
initial condition
$$
\hskip -2em
\ddot\varphi_i\,\hbox{\vrule height 8pt depth 8pt width 0.5pt}_{\,t=0}=0
\tag16.2
$$
in addition to initial conditions \thetag{7.5}. Let's calculate second
derivatives $\ddot\varphi_i$ for deviation functions $\varphi_1,\,\ldots,
\,\varphi_{n-1}$ by differentiating formula \thetag{7.10}. First of all
note that formula \thetag{7.10} itself can be written in terms of
covariant derivatives
$$
\hskip -2em
\dot\varphi_i=-\sum^n_{s=1}\frac{\tilde\nabla^s\!H}{\Omega}\,
\nabla_{\!\boldsymbol\tau_i}p_s-\sum^n_{s=1}\left(\frac{\nabla_{\!s}H}
{\Omega}-Q_s\!\right)\tau^s_i.
\tag16.3
$$
Note that the equations of Newtonian dynamics in form \thetag{7.5} also
admit covariant derivatives instead of partial derivatives in them:
$$
\xalignat 2
&\hskip -2em
\dot x^s=\frac{\tilde\nabla^s\!H}{\Omega},
&&\nabla_{\!t}p_s=-\frac{\nabla_{\!s}H}{\Omega}+Q_s.
\tag16.4
\endxalignat
$$
Applying covariant derivative $\nabla_{\!\boldsymbol\tau_i}$ to
\thetag{16.4} and denoting $\nabla_{\!\boldsymbol\tau_i}\bold p=
\boldsymbol\xi_i$, we obtain the following differential equations
for components of vector $\boldsymbol\tau_i$ and covector
$\boldsymbol\xi_i$:
$$
\align
&\hskip -2em
\nabla_{\!t}\tau^s_i=\sum^n_{r=1}\tilde\nabla^r\!\!
\left(\frac{\tilde\nabla^s\!H}{\Omega}\right)\xi_{ri}+
\sum^n_{r=1}\nabla_{\!r}\!\!\left(\frac{\tilde\nabla^s\!H}
{\Omega}\right)\tau^r_i,
\tag16.5\\
\vspace{2ex}
&\hskip -2em
\aligned
\nabla_{\!t}&\xi_{si}
+\sum^n_{r=1}\sum^n_{q=1}\sum^n_{m=1}D^{mq}_{rs}\,p_m
\left(\frac{\nabla_{\!q}H}{\Omega}-Q_q\!\right)\tau^r_i\,+\\
\vspace{1ex}
&+\sum^n_{q=1}\sum^n_{m=1}\frac{\tilde\nabla^q\!H\,p_m}{\Omega}
\left(\,\shave{\sum^n_{r=1}}\tilde R^m_{sqr}\,\,\tau^r_i
+\shave{\sum^n_{r=1}}D^{mr}_{qs}\,\xi_{ri}\!\right)=\\
\vspace{1ex}
&=-\sum^n_{r=1}\tilde\nabla^r\!\!\left(\frac{\nabla_{\!s}H}
{\Omega}-Q_s\!\right)\xi_{ri}-
\sum^n_{r=1}\nabla_{\!r}\!\!\left(\frac{\nabla_{\!s}H}{\Omega}
-Q_s\!\right)\tau^r_i,
\endaligned
\tag16.6
\endalign
$$
Here $D^{mq}_{rs}$ and $D^{mr}_{qs}$ are components of dynamic
curvature tensor $\bold D$ given by formula
$$
\hskip -2em
D^{kr}_{ij}=-\frac{\partial\Gamma^k_{ij}}{\partial p_r}.
\tag16.7
$$
Tensor $\bold D$ has no analogs in Riemannian geometry since its
components \thetag{16.7} do vanish for non-extended connections.
In \thetag{16.6} we have quantities $\tilde R^m_{sqr}$ which are
components of another curvature tensor $\tilde\bold R$ given by
formula
$$
\hskip -2em
\gathered
\tilde R^k_{rij}=\frac{\partial\Gamma^k_{jr}}{\partial x^i}-
\frac{\partial\Gamma^k_{ir}}{\partial x^j}+\sum^n_{m=1}
\Gamma^k_{im}\,\Gamma^m_{jr}-\sum^n_{m=1}\Gamma^k_{jm}\,
\Gamma^m_{ir}\,+\\
\vspace{1ex}
+\sum^n_{m=1}\sum^n_{\alpha=1}p_\alpha\,\Gamma^\alpha_{mi}\,
\frac{\partial\Gamma^k_{jr}}{\partial p^m}-\sum^n_{m=1}
\sum^n_{\alpha=1}p_\alpha\,\Gamma^\alpha_{mj}\,\frac{\partial
\Gamma^k_{ir}}{\partial p^m}.
\endgathered
\tag16.8
$$
In Riemannian geometry \thetag{16.8} reduces to standard formula
for components of Riemann curvature tensor.\par
    Now let's differentiate \thetag{16.3} with respect to time
variable $t$. It is equivalent to applying covariant derivative
$\nabla_{\!t}$ to to this equality:
$$
\hskip -2em
\gathered
\ddot\varphi_i=-\sum^n_{s=1}\nabla_{\!t}\!\left(
\frac{\tilde\nabla^s\!H}{\Omega}\right)\,\xi_{si}
-\sum^n_{s=1}\frac{\tilde\nabla^s\!H}{\Omega}\,
\nabla_{\!t}\xi_{si}\,-\\
\vspace{1ex}
-\sum^n_{s=1}\nabla_{\!t}\!\left(\frac{\nabla_{\!s}H}{\Omega}
-Q_s\!\right)\tau^s_i-\sum^n_{s=1}\left(\frac{\nabla_{\!s}H}
{\Omega}-Q_s\!\right)
\nabla_{\!t}\tau^s_i.
\endgathered
\tag16.9
$$
Substituting \thetag{16.5} and \thetag{16.6} into \thetag{16.9},
we obtain the following equality for $\ddot\varphi_i$:
$$
\allowdisplaybreaks
\gather
\ddot\varphi_i=\sum^n_{r=1}\left(\,\shave{\sum^n_{s=1}}
\frac{\nabla_{\!s}\Omega}{\Omega^2}\,\frac{\tilde\nabla^s\!H}
{\Omega}+\shave{\sum^n_{s=1}\frac{\tilde\nabla^s\Omega}
{\Omega^2}\left(-\frac{\nabla_{\!s}H}{\Omega}+Q_s\!\right)}
\right)\tilde\nabla^r\!H\,\,\xi_{ri}\,-\\
\vspace{1ex}
-\sum^n_{r=1}\sum^n_{s=1}\frac{\tilde\nabla^s\!H}{\Omega}
\left(\tilde\nabla^rQ_s+\frac{\tilde\nabla^r\Omega}{\Omega}
\,Q_s\!\right)\xi_{ri}+\sum^n_{r=1}\left(\,\shave{\sum^n_{s=1}}
\frac{\nabla_{\!s}\Omega}{\Omega^2}\,\,\frac{\tilde\nabla^s\!H}
{\Omega}\,+\right.\\
\vspace{1ex}
\left.+\shave{\sum^n_{s=1}}\frac{\tilde\nabla^s\Omega}
{\Omega^2}\left(-\frac{\nabla_{\!s}H}{\Omega}+Q_s\!\right)
\right)\nabla_{\!r}H\,\,\tau^r_i+\sum^n_{r=1}\sum^n_{s=1}
\left(\vphantom{\frac{\nabla_{\!r}\Omega\,Q_s}{\Omega}}
\nabla_{\!s}Q_r-\nabla_{\!r}Q_s\,-\right.\\
\vspace{1ex}
\left.-\,\frac{\nabla_{\!r}\Omega}{\Omega}\,Q_s\!\right)
\frac{\tilde\nabla^s\!H}{\Omega}\,\tau^r_i+\sum^n_{r=1}
\sum^n_{s=1}\left(-\frac{\nabla_{\!s}H}{\Omega}+Q_s\!\right)
\tilde\nabla^sQ_r\,\tau^r_i
\endgather
$$
Terms with curvature tensors are canceled due to the following
identities:
$$
\align
&\hskip -2em
[\nabla_{\!i},\,\nabla_{\!j}]H=-\sum^n_{k=1}\sum^n_{s=1}p_k\,
\tilde R^k_{sij}\,\tilde\nabla^s\!H,
\tag16.10\\
\vspace{1ex}
&\hskip -2em
[\nabla_{\!i},\,\tilde\nabla^j]H=\sum^n_{k=1}\sum^n_{s=1}p_k\,
D^{kj}_{is}\,\tilde\nabla^s\!H.
\tag16.11
\endalign
$$
The identities similar to \thetag{16.10} and \thetag{16.11}
in $\bold v$-representation were derived in Chapter~\uppercase
\expandafter{\romannumeral 3} of thesis \cite{6}.\par
    Thus, formula for $\ddot\varphi_i$ is derived (see above).
It is rather huge. In order to simplify this formula we introduce
two extended fields $\boldsymbol\alpha$ and $\boldsymbol\beta$
with components
$$
\gather
\hskip -2em
\gathered
\alpha^r=\left(\,\shave{\sum^n_{s=1}}\frac{\nabla_{\!s}\Omega}
{\Omega^2}\,\frac{\tilde\nabla^s\!H}{\Omega}+\shave{\sum^n_{s=1}
\frac{\tilde\nabla^s\Omega}{\Omega^2}\left(-\frac{\nabla_{\!s}H}
{\Omega}+Q_s\!\right)}\right)\tilde\nabla^r\!H\,-\\
\vspace{1ex}
-\sum^n_{s=1}\frac{\tilde\nabla^s\!H}{\Omega}\left(\tilde\nabla^rQ_s
+\frac{\tilde\nabla^r\Omega}{\Omega}\,Q_s\!\right)
\endgathered
\tag16.12\\
\vspace{2ex}
\hskip -2em
\gathered
\beta_r=\left(\,\shave{\sum^n_{s=1}}
\frac{\nabla_{\!s}\Omega}{\Omega^2}\,\,\frac{\tilde\nabla^s\!H}
{\Omega}+\shave{\sum^n_{s=1}}\frac{\tilde\nabla^s\Omega}
{\Omega^2}\left(-\frac{\nabla_{\!s}H}{\Omega}+Q_s\!\right)
\right)\nabla_{\!r}H\,+\\
\vspace{1ex}
+\sum^n_{s=1}\left(\vphantom{\frac{\nabla_{\!r}
\Omega\,Q_s}{\Omega}}\nabla_{\!s}Q_r-\nabla_{\!r}Q_s
-\,\frac{\nabla_{\!r}\Omega}{\Omega}\,Q_s\!\right)
\frac{\tilde\nabla^s\!H}{\Omega}-\sum^n_{s=1}
\left(\frac{\nabla_{\!s}H}{\Omega}-Q_s\!\right)
\tilde\nabla^sQ_r.
\endgathered
\tag16.13
\endgather
$$
Using notations \thetag{16.12} and \thetag{16.13}, we can write
formula for $\ddot\varphi_i$ as follows:
$$
\hskip -2em
\ddot\varphi_i=\sum^n_{s=1}\alpha^s\,\nabla_{\!\boldsymbol
\tau_i}p_s+\sum^n_{s=1}\beta_s\,\tau^s_i.
\tag16.14
$$
Then, using projector $\bold P$ with components \thetag{12.3}, we can
transform this expression:
$$
\hskip -2em
\gathered
\ddot\varphi_i=\sum^n_{r=1}\sum^n_{s=1}\alpha^r\,P^s_r\,
\nabla_{\!\boldsymbol\tau_i}p_s+\sum^n_{r=1}\sum^n_{s=1}
\beta_r\,P^r_s\,\tau^s_i\,+\\
\vspace{1ex}
+\sum^n_{r=1}\sum^n_{s=1}\alpha^r\,p_r\,
\frac{\tilde\nabla^s\!H}{\Omega}\,\nabla_{\!\boldsymbol\tau_i}p_s
+\sum^n_{r=1}\sum^n_{s=1}\beta_r\,\frac{\tilde\nabla^r\!H}{\Omega}
\,p_s\,\tau^s_i.
\endgathered
\tag16.15
$$
Comparing \thetag{16.15} with formula \thetag{16.3} for $\dot\varphi_i$,
we can write formula \thetag{16.15} as
$$
\gathered
\ddot\varphi_i+\left<\bold p\,|\,\boldsymbol\alpha\right>
\dot\varphi_i
=\sum^n_{r=1}\sum^n_{s=1}\alpha^r\,P^s_r\,
\nabla_{\!\boldsymbol\tau_i}p_s+\sum^n_{r=1}\sum^n_{s=1}
\beta_r\,P^r_s\,\tau^s_i\,-\\
\vspace{1ex}
-\sum^n_{s=1}\left<\bold p\,|\,\boldsymbol\alpha\right>
\left(\frac{\nabla_{\!s}H}{\Omega}-Q_s\!\right)\tau^s_i
+\sum^n_{r=1}\sum^n_{s=1}\beta_r\,\frac{\tilde\nabla^r\!H}{\Omega}
\,p_s\,\tau^s_i.
\endgathered
\tag16.16
$$
Now let's introduce other two extended fields $\sigma$ and
$\boldsymbol\eta$ with components
$$
\xalignat 2
&\hskip -2em
\eta_r=\beta_r-\left<\bold p\,|\,\boldsymbol\alpha\right>
\left(\frac{\nabla_{\!r}H}{\Omega}-Q_r\!\right),
&&\sigma=\sum^n_{r=1}\frac{\tilde\nabla^r\!H}{\Omega}\,
\eta_r.\quad
\tag16.17
\endxalignat
$$
Then, taking into account formula \thetag{2.6} for $\varphi_i$,
we can write \thetag{16.16} as follows:
$$
\ddot\varphi_i+\left<\bold p\,|\,\boldsymbol\alpha\right>
\dot\varphi_i-\sigma\,\varphi_i
=\sum^n_{r=1}\sum^n_{s=1}\alpha^r\,P^s_r\,
\nabla_{\!\boldsymbol\tau_i}p_s+\sum^n_{r=1}\sum^n_{s=1}
\eta_r\,P^r_s\,\tau^s_i.
\tag16.18
$$
Combining \thetag{16.18} with \thetag{7.2} and \thetag{16.2},
we obtain the equality
$$
\hskip -2em
\sum^n_{r=1}\sum^n_{s=1}\alpha^r\,P^s_r\,
\nabla_{\!\boldsymbol\tau_i}p_s+\sum^n_{r=1}\sum^n_{s=1}
\eta_r\,P^r_s\,\tau^s_i=0.
\tag16.19
$$
If we remember operator $\bold b$ determined by \thetag{13.1} and
\thetag{13.2} and if we use formulas \thetag{13.3} and \thetag{13.4}
for bilinear form $b$, then \thetag{16.19} is written as
$$
\hskip -2em
-\left<\bold b(\boldsymbol\tau_i)\,|\,\bold P\boldsymbol\alpha\right>
+\left<\boldsymbol\eta\,|\,\bold P\boldsymbol\tau_i\right>=0.
\tag16.20
$$
Now recall that in the beginning of this section we have taken
a point $q_0=(p_0,\bold p_0)$ and considered a set of hypersurfaces
in $M$ passing through the point $p_0$ tangent to null-space of
covector $\bold p_0\neq 0$. Therefore, when equality \thetag{16.20}
is written for the point $p_0$, we can replace $\bold P\boldsymbol
\tau_i$ by $\bold P\bold X$, where $\bold X$ is an arbitrary vector of
tangent space $T_{p_0}(M)$. In a similar way, due to theorem~13.2,
covector $\bold b(\boldsymbol\tau_i)$ in \thetag{16.20} can be replaced
by arbitrary covector $\bold y\in T^*_{p_0}(M)$. Hence \thetag{16.20}
breaks into two parts
$$
\xalignat 2
&\hskip -2em
\left<\bold y\,|\,\bold P\boldsymbol\alpha\right>=0,
&&\left<\boldsymbol\eta\,|\,\bold P\bold X\right>=0
\tag16.21
\endxalignat
$$
with arbitrary vector $\bold X$ and arbitrary covector $\bold y$.
In coordinate form these two equalities \thetag{16.21} are equivalent
to the following ones:
$$
\xalignat 2
&\hskip -2em
\sum^n_{s=1}P^r_s\,\alpha^s=0,
&&\sum^n_{s=1}P^s_r\,\eta_s=0.
\tag16.22
\endxalignat
$$
Looking at \thetag{16.12}, \thetag{16.13}, and \thetag{16.17}, we
see that the equalities \thetag{16.22} form a system of partial
differential equations for components of covector $\bold Q$ written
in terms of covariant derivatives $\nabla$ and $\tilde\nabla$.
They are called {\it weak normality equations}. Let's write them
explicitly. For the first equation \thetag{16.22} we have
$$
\hskip -2em
\sum^n_{r=1}\sum^n_{s=1}\frac{\tilde\nabla^s\!H}{\Omega}
\left(\tilde\nabla^rQ_s+\frac{\tilde\nabla^r\Omega}{\Omega}
\,Q_s\!\right)P^q_r=0.
\tag16.23
$$
Second equality \thetag{16.22} leads to more huge equations
$$
\hskip -2em
\gathered
\sum^n_{r=1}\sum^n_{s=1}\left(
\left(\nabla_{\!s}Q_r+\frac{\nabla_{\!s}\Omega}{\Omega}\,Q_r
-\nabla_{\!r}Q_s+\frac{\nabla_{\!r}\Omega}{\Omega}\,Q_s\!\right)
\frac{\tilde\nabla^s\!H}{\Omega}\,+\right.\\
\vspace{1ex}
+\sum^n_{m=1}\left(\frac{\nabla_{\!r}H}{\Omega}-Q_r\!\right)\!
\left(\tilde\nabla^mQ_s+\frac{\tilde\nabla^m\Omega}{\Omega}
\,Q_s\!\right)\frac{\tilde\nabla^s\!H}{\Omega}\,p_m\,-
\vspace{1ex}
\left.-\left(\frac{\nabla_{\!s}H}{\Omega}-Q_s\!\right)\!
\left(\tilde\nabla^sQ_r+\frac{\tilde\nabla^s\Omega}{\Omega}
\,Q_r\!\right)\right)P^r_q=0.
\endgathered
\tag16.24
$$
\proclaim{Theorem 16.2} \ Strong normality condition for Newtonian
dynamical system \thetag{16.3} implies weak normality equations
\thetag{16.23} and \thetag{16.24} to be fulfilled at all points
$q=(p,\bold p)$ of cotangent bundle $T^*\!M$, where $\bold p\neq 0$.
\endproclaim
    Note that theorem~16.2 is stated for Newtonian dynamical system
written in form \thetag{16.4}. However, weak normality equations
\thetag{16.23} and \thetag{16.24} are invariant under the transformations
\thetag{15.1} (this can be proved by direct calculations). Therefore
theorem~16.2 is equally applicable to Newtonian dynamical system written
in form of the equations \thetag{7.5}.
\head
17. Equivalence of strong and complete normality conditions.
\endhead
    In section~16 we have derived weak normality equations \thetag{16.23}
and \thetag{16.24} from strong normality condition for Newtonian dynamical
system (see definition~8.1). Here we reveal their relation to weak normality
condition considered in section~6. As in section~6, let's consider
one-parametric family of trajectories $p=p(t,y)$ of Newtonian dynamical
system \thetag{16.4}. In local coordinates it is represented by functions
\thetag{6.1}. Differentiating them with respect to parameter $y$, we obtain
variation vector $\boldsymbol\tau$, see formula \thetag{6.3}. Then we can
define deviation function \thetag{6.4}. Unlike section~6, in the above
calculations in section~16 we used $\bold p$-representation rather than
$\bold v$-representation. Indeed, the equations of Newtonian dynamics
\thetag{16.4} and all normality equations \thetag{14.2}, \thetag{14.12},
\thetag{16.23}, and \thetag{16.24} are written in terms of Hamilton function
$H$ and in terms of covariant derivatives \thetag{10.3} and \thetag{10.7},
which require momentum representation of extended tensor fields. For this
reason, instead of functions $\theta^i(t,y)$ in \thetag{6.5}, we consider
components of vector $\boldsymbol\xi=\nabla_{\!\boldsymbol\tau}\bold p$:
$$
\xi_s=\nabla_{\!\boldsymbol\tau}p_s=\frac{\partial p_s}{\partial y}-
\sum^n_{k=1}\sum^n_{q=1}\Gamma^k_{sq}\,p_k\,\tau^q
$$
(compare with formula \thetag{12.9} above). Functions $\tau^1,\,\ldots,\,
\tau^n,\,\xi^1,\,\ldots,\,\xi^n$ considered as functions of time variable
$t$ for fixed $y$ satisfy a system ordinary differential equations which
are quite the same as the equations \thetag{16.5} and \thetag{16.6} above,
except for the absence of index $i$ now:
$$
\align
&\hskip -2em
\nabla_{\!t}\tau^s=\sum^n_{r=1}\tilde\nabla^r\!\!
\left(\frac{\tilde\nabla^s\!H}{\Omega}\right)\xi_{ri}+
\sum^n_{r=1}\nabla_{\!r}\!\!\left(\frac{\tilde\nabla^s\!H}
{\Omega}\right)\tau^r,
\tag17.1\\
\vspace{2ex}
&\hskip -2em
\aligned
\nabla_{\!t}&\xi_s
+\sum^n_{r=1}\sum^n_{q=1}\sum^n_{m=1}D^{mq}_{rs}\,p_m
\left(\frac{\nabla_{\!q}H}{\Omega}-Q_q\!\right)\tau^r\,+\\
\vspace{1ex}
&+\sum^n_{q=1}\sum^n_{m=1}\frac{\tilde\nabla^q\!H\,p_m}{\Omega}
\left(\,\shave{\sum^n_{r=1}}\tilde R^m_{sqr}\,\,\tau^r
+\shave{\sum^n_{r=1}}D^{mr}_{qs}\,\xi_r\!\right)=\\
\vspace{1ex}
&=-\sum^n_{r=1}\tilde\nabla^r\!\!\left(\frac{\nabla_{\!s}H}
{\Omega}-Q_s\!\right)\xi_r-
\sum^n_{r=1}\nabla_{\!r}\!\!\left(\frac{\nabla_{\!s}H}{\Omega}
-Q_s\!\right)\tau^r.
\endaligned
\tag17.2
\endalign
$$
These equations \thetag{17.1} and \thetag{17.2} can be understood
as linearizations of the equations of Newtonian dynamics \thetag{16.4}.
They are linear equations with respect to functions $\tau^1,\,\ldots,
\,\tau^n,\,\xi^1,\,\ldots,\,\xi^n$, but with non-constant coefficients.
Their coefficients are functions of time variable $t$ determined
by a trajectory $p=p(t)$ of Newtonian dynamical system
\thetag{16.4}. For a fixed trajectory $p=p(t)$, i\.\,e\. when parameter
$y$ in $p=p(t,y)$ is fixed, solutions of the equations \thetag{17.1} and
\thetag{17.2} form $n$-dimensional linear space. We denote this space
by $\goth T$. In essential, it is the same space $\goth T$ as in section~6.
\par
    Let's consider deviation function $\varphi$ determined by formula
\thetag{6.4} and its time derivatives. Now $\varphi$ and functions
$\tau^1,\,\ldots,\,\tau^n,\,\xi^1,\,\ldots,\,\xi^n$ are not related
to any hypersurface $S$. Nevertheless, repeating the same steps as in
deriving formula \thetag{16.18}, one can derive the following equality
for deviation function $\varphi$:
$$
\ddot\varphi+\left<\bold p\,|\,\boldsymbol\alpha\right>
\dot\varphi-\sigma\,\varphi=\sum^n_{r=1}\sum^n_{s=1}\alpha^r\,P^s_r\,
\xi_s+\sum^n_{r=1}\sum^n_{s=1}\eta_r\,P^r_s\,\tau^s.
\tag17.3
$$
Remember that weak normality equations \thetag{16.23} and \thetag{16.24}
are expanded form of the equations \thetag{16.22}. Therefore, if weak
normality equations are fulfilled, then \thetag{17.3} reduces to second
order ordinary differential equation for $\varphi$:
$$
\ddot\varphi+\left<\bold p\,|\,\boldsymbol\alpha\right>
\dot\varphi-\sigma\,\varphi=0.
\tag17.4
$$
Comparing \thetag{17.4} with \thetag{6.11}, we find that Newtonian
dynamical system fits the definition~6.1 under the condition that
weak normality equations \thetag{16.23} and \thetag{16.24} for
extended covector field $\bold Q$ are fulfilled.\par
    Converse result is also valid, i\.\,e\. weak normality condition
stated in definition~6.1 implies weak normality equations \thetag{16.23}
and \thetag{16.24} to be fulfilled. Let's prove it. Suppose that we
take some trajectory of $p=p(t)$ of Newtonian dynamical system
\thetag{16.4}. This determines coefficients in linear differential
equations \thetag{17.1} and \thetag{17.2} for components of
$\boldsymbol\tau$ and $\boldsymbol\xi$ and fixes linear space $\goth T$.
Then, according to definition~6.1, for each solution of the equations
\thetag{17.1} and \thetag{17.2} corresponding deviation function
$\varphi$ should satisfy differential equation \thetag{6.11}. 
For time
derivatives $\dot\varphi$ and $\ddot\varphi$ we have the following
equalities:
$$
\align
&\hskip -2em
\dot\varphi=-\sum^n_{s=1}\frac{\tilde\nabla^s\!H}{\Omega}\,
\xi^s-\sum^n_{s=1}\left(\frac{\nabla_{\!s}H}
{\Omega}-Q_s\!\right)\tau^s,
\tag17.5\\
\vspace{1ex}
&\hskip -2em
\ddot\varphi=\sum^n_{s=1}\alpha^s\,\xi^s
+\sum^n_{s=1}\beta_s\,\tau^s.
\tag17.6
\endalign
$$
(compare with \thetag{16.3} and \thetag{16.14} above). For the
function $\varphi$ itself we have formula \thetag{6.4}. Due to
differential equation \thetag{6.11} from \thetag{17.5} and
\thetag{17.6} we derive
$$
\sum^n_{s=1}\left(\alpha^s+\Cal A\,\frac{\tilde\nabla^s\!H}{\Omega}
\right)\xi^s
+\sum^n_{s=1}\left(\beta_s+\Cal A\left(\frac{\nabla_{\!s}H}
{\Omega}-Q_s\!\right)-\Cal B\,p_s\!\right)\tau^s=0.
$$
According to definition~6.1, this equality should be fulfilled for
all solutions of differential equations \thetag{17.1} and \thetag{17.2}.
Therefore
$$
\align
&\hskip -2em
\alpha^s+\Cal A\,\frac{\tilde\nabla^s\!H}{\Omega}=0,
\tag17.7\\
\vspace{1ex}
&\hskip -2em
\beta_s+\Cal A\left(\frac{\nabla_{\!s}H}
{\Omega}-Q_s\!\right)-\Cal B\,p_s=0.
\tag17.8
\endalign
$$
Multiplying \thetag{17.7} by $p_s$ and summing over $s$ running
from $1$ to $n$, we obtain 
$$
\Cal A=-\sum^n_{s=1}\alpha^s\,p_s=-\left<\bold p\,|\,
\boldsymbol\alpha\right>.
\tag17.9
$$
In a similar way, multiplying \thetag{17.8} by $\tilde\nabla^s\!H$
and summing over $s$ running from $1$ to $n$, we derive formula
for coefficient $\Cal B$ in \thetag{6.11}:
$$
\hskip -2em
\Cal B=\sum^n_{s=1}\frac{\tilde\nabla^s\!H}{\Omega}
\left(\beta_s+\Cal A\left(\frac{\nabla_{\!s}H}{\Omega}
-Q_s\!\right)\right)=\sigma.
\tag17.10
$$
Now, let's substitute \thetag{17.10} and \thetag{17.9} back into
the equations \thetag{17.7} and \thetag{17.8}. Then let's multiply
\thetag{17.7} and \thetag{17.8} by $P^q_s$ and $P^s_q$ respectively
and sum them over index $s$. This yields the equalities \thetag{16.22}
which are equivalent to weak normality equations \thetag{16.23} and
\thetag{16.24}. Thus we have proved the following theorem.
\proclaim{Theorem 17.1} Weak normality condition stated in
definition~6.1, when applied to Newtonian dynamical system written
as \thetag{7.5}, is equivalent to the system of weak normality
equations \thetag{16.23} and \thetag{16.24} that should be fulfilled
at all points $q=(p,\bold p)$ of cotangent bundle $T^*\!M$,
where $\bold p\neq 0$.
\endproclaim
Combining theorems~16.1, 16.2 and 17.1, we obtain another theorem.
\proclaim{Theorem 17.2} Strong and complete normality conditions
are equivalent to each other either for $n=2$ and in multidimensional
case for $n\geqslant 3$.
\endproclaim
\head
18. Summary and conclusions.
\endhead
    Primary goal of present paper is to generalize theory of Newtonian
dynamical systems admitting normal shift from Riemannian geometry to the
geometry determined by some Lagrangian or, equivalently, by some Hamiltonian
dynamical system. In the above sections this goal is reached in essential:
\roster
\widestnumber\item{--}
\rosteritemwd=0pt
\item"--" we have found proper statement for the concept of normal
shift in Lagrangian geometry and have defined class of Newtonian
dynamical systems admitting normal shift of hypersurfaces;
\item"--" we have derived complete system of normality equations
thus obtaining effective tool for studying this class of dynamical
systems.
\endroster
However, some details of constructed theory appear to be different from
those one could expect in the beginning. Indeed, being generalization
of Riemannian geometry, geometry of Lagrangian dynamical system could
have some connection canonically associated with Lagrange function.
But even if such connection does exist, theory of normal shift does not
reveal it. All normality equations \thetag{14.2}, \thetag{14.12},
\thetag{16.23}, and \thetag{16.24} are invariant under transformations
\thetag{15.1}. Therefore they are of connection-free nature. There is
a problem of writing them in coordinate covariant tensorial form without
use of connection at all. This problem will be considered in separate
paper.
\head
19. Acknowledgements.
\endhead
     This work is supported by grant from Russian Fund for Basic
Research (coordinator of project Ya\.~T.~Sultanaev), and by grant
from Academy of Sciences of the Republic Bashkortostan (coordinator
N.~M.~Asadullin). I am grateful to these organizations for financial
support.\par

\enddocument
\end